\documentclass[12pt]{article}
\usepackage[mathscr]{eucal}
\usepackage{amsbsy}
\usepackage{amssymb}
\usepackage{amsmath}
\usepackage{amsthm}
\usepackage{graphics}
\usepackage{epsfig}

\usepackage{latexsym}

\usepackage[left=2cm, top=2cm, right=1.5cm, bottom=2cm] {geometry}

\def\phi{\varphi}

\def\bbr{{\mathbb R}}

\def\bz{{\bf Z}}

\def\bz{{\mathbb Z}} %Keistas Z apribrezimas

\def\Bl{{\bf L}}

\def\bk{{\bf k}}
\def\bt{{\bf t}}

\def\bu{{\bf u}}
\def\bl{{\bf l}}

\def\bn{{\bf n}}

\def\bi{{\bf i}}
\def\Bl1{{\bf 1}}
\def\B2{{\bf 2}}
\def\B0{{\bf 0}}

\def\bm{{\bf m}}

 \newcommand{\NN}{\mathbb{N} }
  \newcommand{\RR}{\mathbb{R} }
  
  \newcommand{\ZZ}{\mathbb{Z} }
 \newcommand{\EE}{\mathbb{E} }
 
\newcommand{\sv}[1]{\lfloor #1 \rfloor}

 \usepackage{color,enumitem}
 \usepackage{multirow,longtable}

\def\al{\alpha}
\def\b{\beta}
\def\d{\delta}

\def\r{\rho}
\def\e{\varepsilon}
\def\g{\gamma}

\def\s{\sigma}
\def\t{\tau}

\def\cl{{\mathscr L}}

\def\=A8{\"o}

\def\fdd{\stackrel{f.d.d.}{\longrightarrow}}

\newcommand{\dd}{{\rm d}}

\newcommand{\ind}[1]{\mathbbm{1}_{#1}}

\newcommand{\beq}{\begin{equation*}}
\newcommand{\eeq}{\end{equation*}}
\newcommand\beqn{\begin{displaymath}}  % no number
\newcommand\eeqn{\end{displaymath}}

\newcommand{\skl}[1]{\left(#1\right)}
 \newcommand{\ri}{{\rm i}}

\newcommand{\abs}[1]{\left|#1\right|}
\newcommand{\brc}[1]{\left\lbrace#1\right\rbrace}

\newcommand{\halmos}{\vspace{3mm} \hfill \mbox{$\Box$}\\[2mm]}

\theoremstyle{plain}
\newtheorem{teo}{Theorem}

\newtheorem{prop}[teo]{Proposition}
\theoremstyle{definition}
\newtheorem{definition}[teo]{Definition}

\usepackage{bbm,tikz,subcaption}
\begin{document}

\title{Some remarks on scaling transition in limit theorems for random fields
\footnotemark[0]\footnotetext[0]{ \textit{Short title:}
Some remarks on scaling transition }
\footnotemark[0]\footnotetext[0]{%
\textit{MSC 2000 subject classifications}. Primary 60G60 .} \footnotemark[0]
\footnotemark[0]\footnotetext[0]{ \textit{Corresponding author:}
Vygantas Paulauskas, Department of Mathematics and
 Informatics, Vilnius university, Naugarduko 24, Vilnius 03225, Lithuania,
 e-mail:vygantas.paulauskas@mif.vu.lt}}

\author{Julius Damarackas$^{\text{\small 1}}$,   Vygantas Paulauskas$^{\text{\small 1}}$ \\
{\small $^{\text{1}}$ Vilnius University, Department of Mathematics
and
 Informatics,}\\}

%\date{}

%\begin{document}

\maketitle

\begin{abstract}

In the paper we present simple examples of linear random fields defined on $\ZZ^2$ and $\ZZ^3$ which exhibit the scaling transition phenomenon. These examples lead to more general definition of the scaling transition and allow to understand the mechanism of appearance of this phenomenon better. In previous papers devoted to the scaling transition it was proved mainly for random fields with finite variance and long-range dependence. We consider random fields with finite and infinite variance. Our examples show that the scaling transition phenomenon can be observed for linear random fields with the so-called negative dependence, which is part of short-range dependence.
Relation of the scaling transition with Lamperti type theorems for random fields is discussed.

\end{abstract}
\vfill
%\eject

\section{Introduction }

Let $\{\xi_i, \ i\ge 1,\}$ be a sequence of random variables and $S_n=\sum_{i=1}^n \xi_i, \ n\ge 1, S_0=0.$ The asymptotic theory of partial sums $S_n$ and related partial sum-processes $S_n(t)=S_{[nt]}, \ t\ge 0$ is well developed and documented in many monographs starting with the classical book \cite{Gnedenko}. Let us note that in the summation theory for sequences we  perform the summation over increasing (with $n\to \infty$) intervals of integers $[1, n]$, although there is  a possibility to consider sums over some increasing sets $A_n\subset A_{n+1}\subset \bz$. Passing from sequences to random fields we have much more flexibility, even in formulation of a problem, and one can say that the asymptotic theory of summation of values of random fields is developed less comparing with the same theory for sequences. To explain our goals and for the simplicity of writing and understanding, we consider a stationary random field (r.f.) $Y=\{Y_{k_1, k_2}, \ (k_1, k_2)\in \bz^2,\}$ indexed by two-dimensional indices. We can perform the summation of the values of the r.f. over some sequence of increasing sets $A_n\subset A_{n+1}\subset \bz^2$, or even over some sets, indexed by multi-indices, and investigate $\sum_{(k_1, k_2)\in A_n}Y_{k_1, k_2}$. Rectangles (sets, indexed by two-dimensional indices) present one of the most simple sets  for such summation, and we can consider
\begin{equation}\label{sum}
S_{n_1, n_2}=S_{n_1, n_2}(Y):=\sum_{k_1=1}^{n_1}\sum_{k_2=1}^{n_2} Y_{k_1, k_2} \ {\rm and} \ S_{n_1, n_2}(t_1, t_2; Y)=S_{[n_1t_1], [n_2t_2]}(Y), \quad t_1\ge 0, \ \ t_2 \ge 0.
\end{equation}
To avoid the centering, let us assume that, if a r.f. $Y$ has the first moment finite, then $EY_{0, 0}=0$. One can expect that under some mild conditions there exist some unboundedly growing constants $A_{n_1, n_2}$ such that finite dimensional distributions (f.d.d.) of $A_{n_1, n_2}^{-1}S_{n_1, n_2}(t_1, t_2)$ converges, as $(n_1, n_2)\to \infty,$ to f.d.d. of some non-trivial (not identically equal to zero) r.f. $V(t, s).$ Here and in what follows $(n_1, n_2)\to \infty$ means that $\min(n_1, n_2)\to \infty.$ It is not difficult to see that the situation described above for r.f. is analogous to the situation for sequences, considered in Lamperti theorems and giving rise to self-similar processes, see \cite{Lamperti}. We refer  a reader to the recent paper \cite{DavPau1} and references therein, where generalizations of Lamperti theorems for r.f. are considered. In our context the following result - Corollary 1 from \cite{DavPau1} - is important. We shall formulate as Proposition only particular case (in \cite{DavPau1} $\bbr^m$-valued random fields on $\bz^d$ are considered).

\begin{prop} \label{prop0} Suppose that $\eta=\{\eta_{i_1, i_2}, (i_1, i_2)\in \bz^2\}$ is a real-valued stationary random field.
If there exists a function $f(n_1, n_2)\to \infty$ as $(n_1, n_2)\to \infty$ and $b\in \bbr$ such that
 \begin{equation}\label{coreq}
\left \{\frac{S_{[n_1t_1], [n_2t_2]}(\eta -b)}{f(n_1, n_2)}, \ (t_1, t_2) \in \bbr_+^2\right \}\fdd \{V(t_1, t_2), \ (t_1, t_2) \in \bbr_+^2 \},
\end{equation}
  where $V$ is a non-degenerate continuous in probability real-valued random field, then $f(a_1, a_2)= a_1^{H_1}a_2^{H_2}L(a_1, a_2)$ with some $H_i>0,  i=1, 2,$ and some coordinate-wise slowly varying function $L,$ and $V$ is $(H_1, H_2)$-multivariate self-similar random field.
\end{prop}
We  do not provide here definitions of a coordinate-wise slowly varying function and $(H_1, H_2)$-multivariate self-similar random field, referring a reader  to \cite{DavPau1}, since for our purposes these notions are unimportant. For us the following fact is important: in the limit theorem, formulated  in Proposition \ref{prop0}, the limit random field does not depend on the way how $(n_1, n_2)$ tends to infinity. Examples of stationary random fields, for which Proposition \ref{prop0} can be applied were given in \cite{Paul20} (see the discussion on directional memory for random fields and scaling transition at the end of this paper). It was shown that if a linear field
\begin{equation}\label{field2}
\eta_{k_1, k_2}=\sum_{i,j=0}^\infty c_{i_1,i_2}\e_{k_1-i_1, k_2-i_2}, \ (k_1, k_2)\in \bz^2,
\end{equation}
where  innovations $\e_{k_1, k_2}, \ (k_1, k_2)\in \bz^2,$ are i.i.d.  random   variables with mean zero and finite variance and the filter  $\{c_{i_1, i_2}, \ i_1\ge 0,  i_2\ge 0\}$ is  of the form $c_{i_1, i_2}=a_{i_1}b_{i_2}$ with some sequences $\{a_i\}, \{b_i\}$, then for such linear  random field Proposition \ref{prop0} can be applied.
On the other hand, there are random fields for which the limit random field for sums $S_{[n_1t_1], [n_2t_2]}(\xi -b)$ depends on the way how $(n_1, n_2)$ tends to infinity.
In the recent papers   \cite{Puplinskaite1}, \cite{Puplinskaite2}, and \cite{Pilipauskaite}    the so-called phenomenon of the scale transition for random fields in the case $d=2$ was described. This phenomenon can be explained  as follows. We use the notation introduced in (\ref{sum}). The sides $n_1$ and $n_2$ of rectangles in (\ref{sum}) are allowed to grow to infinity arbitrary, now let us suppose that these lengths are connected by the relation $n_1=n, \ n_2=n^\g, \ \g>0$.  Let us consider the r.f.

\begin{equation}\label{sumgamma}
Z_{n, \g}(t_1, t_2)=S_{n, n^{\g}}(t_1, t_2), \quad  t_1\ge 0, \quad t_2\ge 0.
\end{equation}
Let us  assume that, for any $\g>0$, there exists a nontrivial random field $V_\g (t_1, t_2)$ and a normalization $A_n (\g) \to \infty$ such that f.d.d. of $A_n^{-1} (\g)Z_{n, \g}(t_1, t_2)$ converges weakly to f.d.d. of $V_\g$. In \cite{Puplinskaite2} the following definition is given.
\begin{definition}\label{sctran}
A random field $Y$ exhibits scaling transition if there exists $\g_0>0$ such that the limit process $V_\g$ is the same, let say $V_+$, for all  $\g >\g_0$ and  another, not obtained by simple scaling, $V_-$, for   $\g <\g_0$. The r.f. $V_{\g_0}$ is called well-balanced scaling limit of $Y$ and r.f. $V_+$ and $V_-$ are called unbalanced  scaling limits of $Y$.
\end{definition}
In \cite{Puplinskaite2} Gaussian random fields are investigated and it is shown that if the spectral density of a Gaussian random field  $Y$ is
\begin{equation}\label{specden}
f(x_1, x_2)=\frac{g(x_1, x_2)}{(|x_1|^2+c|x_2|^{2H_2/H_1})^{H_1/2}}, \quad (x_1, x_2)\in [-\pi, \pi]^2,
\end{equation}
where $0<H_1\le H_2<\infty, \ H_1H_2<H_1+H_2, \ c>0$ and $g$ is bounded and continuous at origin and $g(0, 0)=1,$ then  the random field $Y$ exhibits the scaling transition at $\g_0=H_1/H_2$ and the expressions of unbalanced and well-balanced scaling limits of $Y$ are described. For the precise formulation see Theorem 3.1 in \cite{Puplinskaite2}. In \cite{Puplinskaite1}   the scaling transition phenomenon is demonstrated for aggregated nearest neighbor random coefficients autoregressive r.f. with finite and infinite variance. In \cite{Pilipauskaite} the same phenomenon is shown for non-linear r.f., obtained by taking Appel polynomials of linear r.f. with filter coefficients decaying at possibly different rate in the horizontal and the vertical directions. All these examples are quite complicated and it is difficult to understand what causes the scaling transition phenomenon. In all three above mentioned papers fields exhibiting scaling transition have one main feature - long-range dependence, and one can think that this feature of r.f. is  responsible for the phenomenon of scaling transition. That this is not the case shows the above given example of a linear field (\ref{field2}) with
 the filter  of the form $c_{i_1, i_2}=a_{i_1}b_{i_2}$ with some sequences $\{a_i\}, \{b_i\}$ then the linear  r.f. (\ref{field2}) does not exhibit the scaling transition, despite the fact that it can have long or short range dependence and all sorts of directional memory (positive, zero or negative in any direction). Also in \cite{Puplinskaite2} it is shown that a Gaussian r.f. with long-range dependence and with spectral density of the form
\begin{equation}\label{specden1}
f(x_1, x_2)=\frac{g(x_1, x_2)}{|x_1|^{2d_1}|x_2|^{2d_2}},
\end{equation}
where $0<d_1, d_2 <1/2$ and $g$ has the same properties as in (\ref{specden}), does not exhibit the scaling transition. Thus, long-range dependence, although important in this problem, is not necessary condition for the scaling transition.

The main goal of the present paper is to provide the most simple examples of r.f. exhibiting the scaling transition in order to understand  the mechanism of appearance of the scaling transition better.
We shall consider only linear r.f. (\ref{field2}) and their sums of the form (\ref{sum}), thus in our considerations there are only two factors, which can cause the scaling transition, namely the filter $\{c_{i_1, i_2}\}$, defining dependence structure of the r.f. $X$ and the way how $(n_1, n_2)$ tends to infinity. Let us note that even in the case of linear random fields with finite variance at present we have no complete answer to the question what properties of  filters give us Lamperti type limit theorem, as in Proposition \ref{prop0}, or scaling transition phenomenon.   Here it is worth to mention that in the papers \cite{Puplinskaite1} and \cite{Pilipauskaite} starting point is a linear r.f., but in  \cite{Puplinskaite1}  a linear r.f. with random coefficients of a filter (autoregressive r.f. with random coefficients) is considered and aggregation procedure is applied. Then the limit r.f. is investigated for scaling transition, while in \cite{Pilipauskaite} Appel polynomials of linear r.f. are considered. Therefore, in these examples there are more factors which can be responsible for the scaling transition. Also we demonstrate that our examples can be easily generalized to the case of r.f. indexed by indices in $\bz^d, \ d>2.$ At the final stage of the preparation of the paper we became aware of the publication \cite{Surg}, where scaling transition is investigated for linear random fields on $\bz^3$.  Taking into account the results of this  paper and our simple examples, one can  say that with $d$ increasing the picture of scaling transition becomes more complicated, even for those  simple examples of linear random fields which we consider.

\section{Preliminaries}

Although in the paper we shall investigate mainly the cases of linear r.f. defined on $\ZZ^2$ or on  $\ZZ^3$, in this section we shall derive some formulae and explain main idea in the case of a general linear field defined on $\ZZ^d$. In what follows letters in bold stand for vectors in $\bbr^d$, and inequalities or equalities between them are component-wise. Functions $\vee$, $\sv{\cdot}$, and product of vectors are understood component-wise: $\bn\vee\bm=(n_1\vee m_1,\dots,n_d\vee m_d)$,  $\sv{\bu}=(\sv{u_1},\dots,\sv{u_d})$,  $\bn\bt=(n_1t_1,\dots,n_dt_d)$. Also we shall use the following notation $\sum_{\bk=\B0}^{\sv{\bn \bt}}:=\sum_{k_1=0}^{\sv{n_1 t_1}} \dots \sum_{k_d=0}^{\sv{n_d t_d}}$, here $\B0:=(0,0,\dots,0)$ and dimension of this vector will be clear from context.

We consider a real valued linear r.f. defined on $\ZZ^d$ and the corresponding partial sum r.f.
\begin{equation}\label{field}
X_\bk=\sum_{\bi\in\ZZ_+^d}c_{\bi}\xi_{\bk-\bi}, \quad S_{\bn}(\bt)=\sum_{\bk=\B0}^{\sv{\bn \bt}} X_{\bk}.
\end{equation}
Here $\xi_{\bi}, \ \bi \in \ZZ^d$  are i.i.d. random variables with the characteristic function (ch.f.) $\exp(-|t|^\al), \ 0,\al \le 2$.
We are interested in limit behavior of appropriately normalized r.f. $A_{\bn}^{-1}S_{\bn}(\bt)$, as $\bn \to \infty$, which means $\min(n_i,i=1,\dots,d)\rightarrow\infty$. Firstly, we shall rewrite the expression of $S_{\bn}(\bt)$ in the following way:
\begin{align*}
S_{\bn}(\bt)
	&=\sum_{\B0\leq \bk\leq \bn\bt} X_{\bk}=\sum_{\B0\leq \bk\leq \bn\bt} \sum_{\bi\geq \B0}c_{\bi}\xi_{\bk-\bi}\\
	&=\sum_{\bk\in\ZZ^d}
	\sum_{\bi\in\ZZ^d}c_{\bi}\xi_{\bk-\bi}\ind{[\B0\leq\bk\leq\bn\bt]}\ind{[\bi\geq\B0]}=\sum_{\bk\in\ZZ^d}
	\sum_{\bl\in\ZZ^d}c_{\bk-\bl}\xi_{\bl}\ind{[\B0\leq\bk\leq\bn\bt]}\ind{[-\bl\geq-\bk]}\\
	&=\sum_{\bk\in\ZZ^d}
	\sum_{\bl\in\ZZ^d}c_{\bk-\bl}\xi_{\bl}\ind{[-\bl\leq\bk-\bl\leq\bn\bt-\bl]}\ind{[\bk-\bl\geq      \B0]}=\sum_{\bk\in\ZZ^d}
	\sum_{\bl\in\ZZ^d}c_{\bk}\xi_{\bl}\ind{[-\bl\leq\bk\leq\bn\bt-\bl]}\ind{[\bk\geq      \B0]}\\
	&=\sum_{\bk\in\ZZ^d}\sum_{\bl\in\ZZ^d}c_{\bk}\xi_{\bl}\ind{[(-\bl)\vee \B0\leq\bk\leq\bn\bt-\bl]}=\sum_{\bl\in\ZZ^d}\xi_{\bl}
	\sum_{\bk\in\ZZ^d}c_{\bk}\ind{[(-\bl)\vee \B0\leq\bk\leq\bn\bt-\bl]}.
\end{align*}
Taking $m$ points $\bt^{(l)}, \ l=1, \dots, m$ we can write  ch.f. of $\sum_{l=1}^{m}x_lS_\bn(\bt^{(l)})$:
\begin{align*}
\EE \exp\left(\ri\sum_{j=1}^{m}x_jS_\bn(\bt^{(j)})\right)&= \EE \exp\left(
\ri\sum_{\bl\in\ZZ^d}\xi_{\bl}\sum_{j=1}^{m}x_j
\sum_{\bk\in\ZZ^d}
c_{\bk}\ind{[(-\bl)\vee \B0\leq\bk\leq\bn\bt^{(j)}-\bl]}\right)\\
&=  \exp\left(-\sum_{\bl\in\ZZ^d}\abs{
\sum_{j=1}^{m}x_j
\sum_{\bk\in\ZZ^d}
c_{\bk}\ind{[(-\bl)\vee \B0\leq\bk\leq\bn\bt^{(j)}-\bl]}}^\al\right).
\end{align*}
We see that we must investigate the quantity
\begin{equation*}
J_{\bn}=J_{\bn}(x_1,\dots,x_m,\bt^{(1)},\dots,\bt^{(m)})=\sum_{\bl\in\ZZ^d}\abs{
	\sum_{j=1}^{m}x_j
	\sum_{\bk\in\ZZ^d}
	c_{\bk}\ind{[(-\bl)\vee \B0\leq\bk\leq\bn\bt^{(j)}-\bl]}}^\al.
\end{equation*}
We write this quantity as integral
\begin{align}\label{generalJn}
J_{\bn}&=\int_{\RR^d}\abs{
	\sum_{j=1}^{m}x_j
	\sum_{\bk\in\ZZ^d}
	c_{\bk}\ind{[(-\sv{\bu})\vee \B0\leq\bk\leq\bn\bt^{(j)}-\sv{\bu}]}}^\al \dd \bu\\ \nonumber
&=\int_{\RR^d}\abs{
	\sum_{j=1}^{m}x_j
	\sum_{\bk\in\ZZ^d}
	c_{\bk}\ind{[(-\sv{\bn\bu})\vee \B0\leq\bk\leq\bn\bt^{(j)}-\sv{\bn\bu}]}}^\al \dd \bn\bu\\ \nonumber
&=\left(\prod_{i=1}^{d}n_i\right)\int_{\RR^d}\abs{
	\sum_{j=1}^{m}x_j
	\sum_{\bk\in\ZZ^d}
	c_{\bk}\ind{[(-\sv{\bn\bu})\vee \B0\leq\bk\leq\bn\bt^{(j)}-\sv{\bn\bu}]}}^\al \dd \bu\\ \nonumber
&=\left(\prod_{i=1}^{d}n_i\right)\int_{\RR^d}h_{\bn}\left(\bu\right)\dd \bu.
\end{align}
Here
\begin{equation}\label{generalhn}
h_{\bn}\left(\bu\right)
 = \abs{\sum_{j=1}^{m}x_j f_{\bn}\left(\bu,\bt^{(j)}\right) }^\al ,
\end{equation}
and
\begin{equation}\label{generalfbn}
f_{\bn}\left(\bu,\bt^{(j)}\right)=\sum_{\bk\in\ZZ^d}
	c_{\bk}\ind{[(-\sv{\bn\bu})\vee \B0\leq\bk\leq\bn\bt^{(j)}-\sv{\bn\bu}]}.
\end{equation}

These three formulae (\ref{generalJn})-(\ref{generalfbn}) will be starting point in the investigation of particular cases $d=2$ or $d=3$ with specific filter coefficients $c_{\bi}$. Main step will be to prove the point-wise convergence of the appropriately normalized  function $h_{\bn}$ to some function $h$ and to show that $h_{\bn}$ is bounded from above by integrable function, then application of Lebesgue dominated convergence theorem will yield the convergence of the normalized quantity $J_{\bn}$.

Since in the examples, which we shall consider to get the scaling transition, the main idea is to take filters with non-zero coefficients  on axes (or other lines), we shall provide the expression of function  (\ref{generalfbn}) in this case. Suppose that we have $d$ sequences $a_q(i), i\in\NN,q=1,\dots,d,$, and let us define
\begin{equation}\label{genci}
c_{\bi}=\sum_{q=1}^{d} a_q(i_q)\ind{[i_q\geq0\text{ and }i_l=0,l\neq q]},
\end{equation}
i.e.,
\begin{equation*}
c_{\bi}=\begin{cases}
a_q(i_q), \text{ if for some}\ q: i_q>0\text{ and } i_l=0,l\neq q,\\
\sum_{q=1}^{d}a_q(0), \text{ if } i_l=\dots =i_d =0,\\
0,\text{ elsewhere}.
\end{cases}
\end{equation*}
Using the notation
\begin{equation*}
A(a,b)=\{ k\in\ZZ: (-a)\vee 0 \leq k \leq b-a \},
\end{equation*}
for this filter we can write the function (\ref{generalfbn}) as follows
\begin{align*}
f_{\bn}\left(\bu,\bt^{(j)}\right)&=\sum_{\bk\in\ZZ^d}
\sum_{q=1}^{d} a_q(k_q)\ind{[k_q\geq0\text{ and }k_l=0,l\neq q]}\ind{[(-\sv{\bn\bu})\vee \B0\leq\bk\leq\bn\bt^{(j)}-\sv{\bn\bu}]}\\
&=\sum_{q=1}^{d} \sum_{\bk\in\ZZ^d}
	 a_q(k_q)
	 \ind{[k_q\geq0\text{ and }k_l=0,l\neq q]}
	 \prod_{i=1}^{d} \ind{A(\sv{n_iu_i},n_it_i^{(j)})}(k_i)\\
	 &=\sum_{q=1}^{d} \sum_{k_q=0}^{\infty}
	 a_q(k_q)
	 \ind{A(\sv{n_qu_q},n_qt_q^{(j)})}(k_q)
      \prod_{i=1,\dots,q-1,q+1,\dots,d}
	 \ind{A(\sv{n_iu_i},n_it_i^{(j)})}(0)\\
&=\sum_{q=1}^{d} \sum_{k_q\in\ZZ}
a_q(k_q)\ind{[(-\sv{n_q u_q})\vee 0\leq k_q\leq n_q t_q^{(j)}-\sv{n_q u_q}]}
\prod_{i=1,\dots,q-1,q+1,\dots,d}
\ind{A(\sv{n_iu_i},n_it_i^{(j)})}(0).
\end{align*}
Denoting
\begin{equation*}
U_q(n_q,u_q,t_q^{(j)})=\sum_{k_q\in\ZZ}
a_q(k_q)\ind{[(-\sv{n_q u_q})\vee 0\leq k_q\leq n_q t_q^{(j)}-\sv{n_q u_q}]}.
\end{equation*}
 we can rewrite the last formula as follows:
\begin{equation}\label{fn2}
f_{\bn}\left(\bu,\bt^{(j)}\right)=\sum_{q=1}^{d}U_q(n_q,u_q,t_q^{(j)})\prod_{i=1,\dots,q-1,q+1,\dots,d}
\ind{A(\sv{n_iu_i},n_it_i^{(j)})}(0).
\end{equation}
Concerning the indicator functions we have the following simple relations:
\begin{equation}\label{indfunc}
\ind{A(\sv{n_iu_i},n_it_i^{(j)})}(0)= \ind{[(-\sv{n_i u_i})\leq 0\leq n_i t_i^{(j)}-\sv{n_i u_i}]}\rightarrow
\ind{[-u_i\le 0\le t_i^{(j)}- u_i]}=\ind{[0,t_i^{(j)}]}(u_i),
\end{equation}
and
\begin{equation}\label{indfunc1}
\ind{A(\sv{n_iu_i},n_it_i^{(j)})}(0)%= \ind{[(-\sv{n_i u_i})\leq 0\leq n_i t_i^{(j)}-\sv{n_i u_i}]}
\leq  \ind{[-n_i u_i \leq 0 \leq n_i t_i^{(j)}-{n_i u_i}+n_i]}
=\ind{[- u_i \leq 0 \leq  t_i^{(j)}-{ u_i}+1]}
=\ind{[ 0,  t_i^{(j)}+1]}(u_i).
\end{equation}
Thus, we see that the main step in investigation of (\ref{fn2}) is to find asymptotic of sums $U_q(n_q,u_q,t_q^{(j)})$, assuming some regular behavior of sequences $a_q(i)$ and some relations between growth of $n_q, \ q=1, \dots , d.$ We shall not continue investigation of the general case $d\ge 2$, since in this way we shall get quite complicated picture, the particular cases $d=2$ and $d=3$ will be more informative and visual. We end this section providing the result about asymptotic behavior of sums, present in (\ref{fn2}). Namely, we assume that we have a sequence $a_i=(1+i)^{-\g}, \ i\ge 1$, and $a_0$ will be defined separately. Here it is appropriate to note, that it is possible to consider more general case $a_i\sim (1+i)^{-\g}L(i),$ where $L$ is a slowly varying function, but since we consider rather specific filters, such generality is unimportant. Let
\begin{equation}\label{Utgamma}
U_{t,\g}(i,n)= \sum_{k=\left( -i \right)\vee 0 }^{\sv{nt}-i}a_{k},
\end{equation}
and our goal is  to find point-wise convergence (for a fixed $t$)
\begin{equation}\label{bound}
\frac{U_{t,\g}(\sv{nu},n)}{z_{\g,n}}\rightarrow H_{\g}(u,t)
\end{equation}
with some sequence $z_{\g,n}$. In order to apply Lebesque dominated convergence theorem, we must get the bound
\begin{equation}\label{bound1}
\frac{U_{t,\g}(\sv{nu},n)}{z_{\g,n}}\leq  G_{\g}(u,t),
\end{equation}
where $G_{\g}(u,t)$ is a function, for a fixed $t$ satisfying
\begin{equation}\label{integrableBound}
\int_{-\infty}^{\infty} \abs{G_{\g}(u,t)}^\al \dd u<\infty.
\end{equation}
Let us define
\begin{equation}\label{sekaNormavimui}
z_{n,\g}=\begin{cases}
1, \text{ if } \g>1 \text{ and }\sum_{j=0}^{\infty}a_j\neq 0,\\
n^{1-\g}, \text{ if } 1<\g<1+1/\al \text{ and }\sum_{j=0}^{\infty}a_j= 0,\\
{\ln n}, \text{ if } \g=1,\\
n^{1-\g}, \text{ if } 1/\al<\g<1.\\
\end{cases}
\end{equation}
Also we define  the function
\begin{equation}\label{Hdef}
H_{\g}(u,t)=\begin{cases}
\sum_{k=0}^{\infty}a_k\ind{[0,t)}(u)\text{, if }\g>1\text{ and }\sum_{j=0}^{\infty}a_j\neq 0,\\
\left ((t-u)_+^{1-\g}-(-u)_+^{1-\g}\right ) (1-\g)^{-1}
\text{, if } 1<\g<1+{1}/{\al} \text{ and }\sum_{j=0}^{\infty}a_j= 0,\\
\ind{[0,t)}(u)\text{, if }\g=1,\\
\left ((t-u)_+^{1-\g}-(-u)_+^{1-\g}\right ) (1-\g)^{-1}
\text{, if }1/\al<\g<1.
\end{cases}
\end{equation}
Here and in what follows we  use the notation $(\cdot)_+=\max(0,\cdot)$.

\begin{prop}\label{prop1} For a sequence  $ \{ a_i, \ i\ge 0,\}$ defined above,  we have the relations (\ref{bound}) - (\ref{integrableBound}) with the functions $z_{n,\g}$ and $H_{\g}(u,t),$ defined in (\ref{sekaNormavimui}) and (\ref{Hdef}), respectively.
Expression of the function  $G_{\g}(u,t)$ is given in (\ref{gfunct}), (\ref{Ggammadef}), and (\ref{gamma1est}).
\end{prop}

{\it Proof of Proposition \ref{prop1}}. %Since in our preprint \cite{PaulDam3} this proposition is proved with all details, here we provide only main steps.
 We start with the case $\g>1$ and $\sum_{k=0}^{ \infty} a_{k}\neq 0$. It is easy to see that for a fixed value of $u$ and a.e.
\begin{equation*}
U_{t,\g}(\sv{nu},n)=\sum_{k=\left( -\sv{nu} \right)\vee 0 }^{\sv{nt}-\sv{nu}}a_{k}
\rightarrow \ind{[0,t)}(u)\sum_{k=0}^{ \infty} a_{k} =:H_{\g}(u,t),
\end{equation*}
we added a.e. since for $u=t$ this limit is $a_0$.
Thus, we have (\ref{bound}) with $z_{\g,n}=1$. For $u\in(-1,t+1)$ we have
\begin{equation*}
\sum_{k=\left( -\sv{nu} \right)\vee 0 }^{\sv{nt}-\sv{nu}}a_{k}
\leq  \ind{(-1,t+1)}(u)\left (|a_0|+ \sum_{k=1}^{ \infty} a_{k}\right ),
\end{equation*}
while for $u\leq -1$ we get
\begin{equation*}
\sum_{k=\left( -\sv{nu} \right)\vee 0 }^{\sv{nt}-\sv{nu}}a_{k}=
\ind{(-\infty,-1]}(u)\int_{-\sv{nu}  }^{\sv{nt}-\sv{nu}+1}(1+{\sv{v}})^{-\g}\dd v\leq
\ind{(-\infty,-1]}(u)\int_{-{nu}  }^{{nt}-{nu}+2n}v^{-\g}\dd v
\end{equation*}
\begin{equation*}
\leq \ind{(-\infty,-1])}(u)\frac{\left( -{u} \right)^{1-\g}-\left( {t}-{u}+2 \right)^{1-\g}}{\g-1}.
\end{equation*}
Denoting
\begin{equation}\label{gfunct}
G_{\g}(u,t)=\begin{cases}
\ind{(-1,t+1)}(u)\left (|a_0|+ \sum_{k=1}^{ \infty} a_{k}\right ),\text{ if }u\in(-1, \infty),\\
\ind{(-\infty,-1]}(u){\left( \left( -{u} \right)^{1-\g}-\left( {t}-{u}+2 \right)^{1-\g} \right)}/{\left( \g-1 \right)},\text{ if }u\leq -1,
\end{cases}
\end{equation}
we get the function, satisfying \eqref{bound1} and \eqref{integrableBound}.

In the case $\g<1$ we have, for  $u\leq t +1$ (otherwise the sum is $0$),
\begin{equation*}
\sum_{k=(-\sv{nu})\vee 0}^{\sv{nt}-\sv{nu}}a_{k}=
\ind{(-\infty,\frac{\sv{nt}}{n}+\frac{1}{n})}(u)\int_{(-\sv{nu})_+}^{\sv{nt}-\sv{nu}+1}a_{\sv{v }} \dd v=
n^{1-\g}\int_{-\infty}^{\infty}\kappa_{\g}(v; u,t) \dd v,
\end{equation*}
where
\begin{equation*}
\kappa_{\g}(v; u,t)=\ind{(-\infty,\frac{\sv{nt}}{n}+\frac{1}{n})}(u) \ind{\left( \frac{(-\sv{nu})_+}{n}, \frac{\sv{nt}-\sv{nu}+1}{n} \right)}(v)
n^{\g}a_{\sv{nv }}.
\end{equation*}
%It is not difficult to see that for $u=t$ we have
%\begin{equation*}
%\kappa_{\g}(v; t,t)= \ind{\left( 0, \frac{1}{n}\right )}(v)n^{\g}a_{\sv{nv }} \to 0, \ {\rm as } \  n\to \infty.
%\end{equation*}
For $u\le t+1$ and  a fixed value of  $v$, we have, as $n\rightarrow\infty$, a.e.
\begin{equation*}
\kappa_{\g}(v; u,t)\rightarrow \ind{(-\infty,t)}(u)\ind{\left( (-u)_+, (t-u) \right)}(v)v^{-\g}.
\end{equation*}
It is easy to see that $\abs{\kappa_{\g}(v; u,t)}\leq \max (1, |a_0|)\ind{(-\infty,t)}(u)\ind{\left( (-u)_+, (t-u+2 ) \right)}(v)v^{-\g}.$
This majorizing function is integrable, therefore, with $z_{\g,n}=n^{1-\g}$,
\begin{equation*}
\frac{U_{t,\g}(\sv{nu},n)}{z_{\g,n}}=n^{\g-1}\sum_{k=(-\sv{nu})_+}^{\sv{nt}-\sv{nu}}a_{k}=\int_{-\infty}^{\infty}
\kappa_{\g}(v; u,t) \dd v\rightarrow
\ind{(-\infty,t)}(u)
\int_{(-u)_+}^{t-u}
v^{-\g} \dd v
\end{equation*}
\begin{equation*}
=\frac{1}{1-\g}\left( ((t-u)_+)^{1-\g}-((-u)_+)^{1-\g} \right)=H_{\g}(u,t).
\end{equation*}
Since
\begin{equation*}
\frac{U_{t,\g}(\sv{nu},n)}{z_{\g,n}}=\int_{-\infty}^{\infty} \kappa_{\g}(v; u,t) \dd v
\leq \max (1, |a_0|)\ind{(-\infty,t+1)}(u) \int_{(-u)_+}^{ t-u+2 } v^{-\g}\dd v,
\end{equation*}
we denote
\begin{equation}\label{Ggammadef}
G_{\g}(u,t):=\max (1, |a_0|)\frac{\ind{(-\infty,t+1)}(u)}{1-\g}\left( (t-u+2)^{1-\g}-((-u)_+)^{1-\g} \right).
\end{equation}
It is easy to note that
\begin{equation*}
\int_{-\infty}^{\infty} \abs{G_{\g}(u,t)}^{\al}\dd u<\infty,
\end{equation*}
thus we have \eqref{bound1} and \eqref{integrableBound}.

The case $\g=1$. Assuming $u\le t+1$ (otherwise the sum is zero), we separate one term:
\begin{equation}\label{sepsum}
\sum_{k=(-\sv{nu})\vee 0}^{\sv{nt}-\sv{nu}}a_{k}=
%a_{(-\sv{nu})\vee 0} + \sum_{k=\left( (-\sv{nu})\vee 0 \right)+1}^{\sv{nt}-\sv{nu}}a_{k}=
a_{(-\sv{nu})\vee 0} + \sum_{k= (1-\sv{nu})\vee 1 }^{\sv{nt}-\sv{nu}}a_{k}.
\end{equation}
For $-2\leq u\le t+1$, we can estimate $a_{(-\sv{nu})\vee 0}\leq C$, while for  $u<-2$ we have
$a_{(-\sv{nu})\vee 0}\leq (-nu)^{-1} \leq (-u)^{-1}$.
Therefore, we have
 \begin{equation}\label{sepsum1}
 (\ln n)^{-1}a_{(-\sv{nu})\vee 0}\rightarrow 0 \quad \text{ and }\quad (\ln n)^{-1}a_{(-\sv{nu})\vee 0} \leq R(u)
 \end{equation}
 with
\begin{equation*}
R(u)=\begin{cases}
C, \text{ if }-2\leq u\le t+1,\\
(-u)^{-1}, \text{ if } u<-2.
\end{cases}
\end{equation*}
 For the separated sum in (\ref{sepsum}), using the change of variables, we can write
\begin{equation*}
\sum_{k=(1-\sv{nu})\vee 1}^{\sv{nt}-\sv{nu}}a_{k}=
\ind{(-\infty,\frac{\sv{nt}}{n})}(u)\int_{(1-\sv{nu})\vee 1}^{\sv{nt}-\sv{nu}+1}a_{\sv{v}}\dd v
=\ln n \int_{-\infty}^{\infty}{\bar \kappa}_{1}(v;u,t)\dd v,
%\ind{(-\infty,t)}(u)\int_{\ln\left( (1-\sv{nu})\vee 1 \right)}^{\ln\left( \sv{nt}-\sv{nu}+1 \right)}a_{\sv{\exp(v)}}\dd \exp(v)
\end{equation*}
where
\begin{equation*}
{\bar \kappa}_{1}(v;u,t)=\ind{(-\infty,\frac{\sv{nt}}{n})}(u)\ind{\left( \frac{\ln\left( (1-\sv{nu})\vee 1 \right)}{\ln n}, \frac{\ln\left( \sv{nt}-\sv{nu}+1 \right)}{\ln n}\right)}(v)
a_{\sv{\exp(v\ln n)}}\exp(v\ln n).
\end{equation*}
We have the point-wise convergence
${\bar \kappa}_{1}(v;u,t)\rightarrow \ind{(0,t)}(u)\ind{\left( 0, 1\right)}(v)$,
and, since $a_{\sv{\exp(v\ln n)}}\exp(v\ln n)\leq 1$, we have the following bound for the integrand
\begin{equation*}\label{gammaEQ1kappaBound}
{\bar \kappa}_{1}(v;u,t) \leq \ind{\left( \frac{\ln\left( (1-\sv{nu})\vee 1 \right)}{\ln n}, \frac{\ln\left( \sv{nt}-\sv{nu}+1 \right)}{\ln n}\right)}(v)\leq
\ind{ \left( 0, 1+ \ln\left( {t}-{u}+2\right)   \right) }(v).
\end{equation*}
Therefore, by the dominated convergence theorem
\begin{equation}\label{sepsum2}
\frac{1}{\ln n}\sum_{k=(1-\sv{nu})\vee 1}^{\sv{nt}-\sv{nu}}a_{k}=
\int_{-\infty}^{\infty} {\bar \kappa}_{1}(v;u,t) \dd v\rightarrow \ind{(0,t)}(u)\int_{0}^{1} 1 \dd v=
\ind{(0,t)}(u)=H_{1}(u,t).
\end{equation}
Applying (\ref{sepsum} we can write
\begin{eqnarray*}
Q_n(u,t) &:=& \frac{1}{\ln n}\sum_{k=(1-\sv{nu})\vee 1}^{\sv{nt}-\sv{nu}}a_{k}\leq
\int_{-\infty}^{\infty}
\ind{\left( \frac{\ln\left( (1-\sv{nu})\vee 1 \right)}{\ln n}, \frac{\ln\left( \sv{nt}-\sv{nu}+1 \right)}{\ln n}\right)}(v)
\dd v \\
   &=& \frac{\ln\left( \sv{nt}-\sv{nu}+1 \right)}{\ln n}-
\frac{\ln\left( (1-\sv{nu})\vee 1 \right)}{\ln n}\\
  &\leq&  \frac{\ln\left( {nt}-{nu}+n \right)-
	\ln\left( (1-\sv{nu})\vee 1 \right)}{\ln n}.
\end{eqnarray*}
For $-2\leq u\le t+1$ we can estimate
\begin{equation*}
Q_n(u,t) \leq \frac{\ln\left( {nt}-{nu}+2n \right) }{\ln n}\leq 1+\ln\left( {t}-{u}+2 \right),
\end{equation*}
while, for $u<-2$, we estimate as follows:
\begin{eqnarray*}
Q_n(u,t) &\leq&  \frac{\ln\left( {nt}-{nu}+n \right)-
	\ln\left( 1-\sv{nu} \right)}{\ln n}
\leq \frac{ \ln\left( {t}-{u}+1 \right)-
	\ln\left( -{u} \right)}{\ln n}\\
 &=&  \ln\left( \frac{{t}-{u}+1 }{-{u}}\right)=
\ln\left( 1+\frac{{t}+1 }{-{u}}\right).
\end{eqnarray*}
It is easy to see that the bounding function decays as $(-u)^{-1}$, for $u\rightarrow-\infty$. Let us denote
\begin{equation}\label{gamma1est}
G_{1}(u,t)=\left(1+\ln\left( {t}-{u}+2 \right)+R(u)\right)\ind{[-2, t+1]}(u)+
\left (\ln\left( 1+\frac{{t}+1 }{-{u}}\right)+R(u)\right )\ind{(-\infty, -2)}(u).
\end{equation}

Since we consider the case $\g=1$, this means that $\al>1$, therefore the function $|G_{1}(u,t)|^\al$ is integrable.
Thus, (\ref{sepsum1}) and (\ref{sepsum2}) gives us (\ref{bound}) with $z_{\g,n}=\ln n$, and we have (\ref{bound1}) and (\ref{integrableBound}) with the function given in (\ref{gamma1est}).

In the case $1<\g<1+{1}/{\al}$ and $\sum_{j=0}^{\infty}a_j= 0$ the proof goes along the same lines as in the case $\g<1$, only we use the equality $\sum_{j=0}^{n}a_j=-\sum_{j=n+1}^{\infty}a_j$. Therefore, we omit the details.
\halmos

\section{The case $d=2$}

In the previous section we derived formulae (\ref{generalJn})-(\ref{generalfbn}), which are the staring point investigating the scaling transition in particular cases $d=2$ and $d=3$. Taking specific filter coefficients $c_{\bi}$ we shall demonstrate what structure and properties of these coefficients exhibit the scaling transition. We recall that we consider r.f. and corresponding partial sum r.f.
\begin{equation}\label{fieldd2}
X_\bk=\sum_{\bi\in\ZZ_+^2}c_{\bi}\xi_{\bk-\bi}, \ \bk=(k_1, k_2)\in \ZZ^2 \ \text {and} \ S_{\bn}(\bt)=\sum_{\bk=\B0}^{\sv{\bn \bt}} X_{\bk}.
\end{equation}

\subsection{Example 1}

Let us take the following two sequences:  $a_j(0)=0, a_j(i)=(1+i)^{-\g_j},  \ i\ge 1, \ \g_j>1/\al$, \ j=1, 2. Consider linear r.f. (\ref{fieldd2})  with the following filter
\begin{equation*}
c_{i_1,i_2}=
\begin{cases}
a_1(i_1), \text{ if } i_1\geq 0, i_2=0,\\
a_2(i_2), \text{ if } i_2\geq 0, i_1
=0,\\
0, \text{ elsewhere}.
\end{cases}
\end{equation*}

This can be written as
\begin{equation}\label{cij}
c_{i_1,i_2}=a_1(i_1)\ind{i_1\geq 0}\ind{i_2=0}+a_2(i_2)\ind{i_2\geq 0}\ind{i_1=0},
\end{equation}
i.e.,  we have  the same filter (\ref{genci}) (with $d=2$) considered in the previous section. From (\ref{generalJn})-(\ref{generalfbn}) and (\ref{fn2}) we have that
\begin{equation}\label{chfuncd2}
\EE \exp\left(\ri ,A_{n_1, n_2}^{-1}\sum_{j=1}^{m}x_jS_\bn(\bt^{(j)})\right)= \exp\left( -A_{n_1,n_2}^{-\al}J_{n_1,n_2}\right),
\end{equation}
where
 \begin{equation}\label{Jbnd2}
J_{n_1,n_2}=n_1n_2\int_{-\infty}^{\infty}\int_{-\infty}^{\infty}f_{n_1,n_2}^{(1)}(u_1,u_2){\dd u}_1{\dd u}_2
\end{equation}
and
\begin{equation}\label{fn1n2}
f_{n_1,n_2}^{(1)}(u_1,u_2)=\abs{\sum_{l=1}^{m}x_l\left( \ind{\{0\leq n_2u_2\leq  \sv{n_2t_2^{(l)}}\}}
	U_{t_1^{(l)},\g_1}(\sv{n_1u_1},n_1)
	+
	\ind{\{0\leq n_1u_1\leq  \sv{n_1t_1^{(l)}}\}}
	U_{t_2^{(l)},\g_2}(\sv{n_2u_2},n_2)  \right) }^\al.
\end{equation}
 We recall (see (\ref{Utgamma})) that
 \begin{equation}\label{stgama}
U_{t_j^{(l)},\g_j}(i,n_j)= \sum_{k=\left( -i \right)\vee 0 }^{\sv{n_jt_j^{(l)}}-i}a_{j}(k), \ j=1, 2.
\end{equation}

Now we can apply Proposition \ref{prop1} and in a standard way (applying Lebesgue dominated convergence theorem) we can get the convergence of ch.f. from (\ref{chfuncd2}).
Since in the expression of the function $f_{n_1,n_2}(u_1,u_2)$ there are two terms, normalization of which can be performed by sequences $z_{\g_1,n_1}$ and $z_{\g_2,n_2}$, see (\ref{sekaNormavimui}), it is clear that it will be impossible to get joint normalization for $f_{n_1,n_2}(u_1,u_2)$ in general case of parameters $\g_1, \g_2$ and  $(n_1,n_2) \to \infty$. But in the case where both normalizing sequences $z_{\g_1,n_1}$ and $z_{\g_2,n_2}$ are equal to $1$, namely, if $\g_i>\max (1, 1/\al), \ i=1, 2$, then we get the following result. Let us denote by $M_\al$ symmetric $\al$-stable measure on $\RR^2$ with Lebesgue control measure and let $\fdd$ stand for the convergence of f.d.d..

\begin{prop}\label{prop2} Suppose that we have a sum (\ref{fieldd2}) of values of a linear r.f. with the filter (\ref{cij}). If $\g_i>\max (1, 1/\al), \ i=1, 2$, then, as $(n_1,n_2) \to \infty$,

\begin{equation}\label{lim1}
(n_1n_2)^{-1/\al}S_{\bn}(\bt)\fdd \sum_{k=1}^\infty (a_1(k)+a_2(k)) \int_0^{t_1}\int_0^{t_2} M_\al (\dd u_1\dd u_2)
\end{equation}
\end{prop}
This proposition means that linear r.f. $\{X_{k_1, k_2}, (k_1, k_2)\in \ZZ^2,  \}$ with such filter does not exhibit the scaling transition and limit in (\ref{lim1}) does not depend on the way how $(n_1, n_2)$ tends to infinity, i.e. the Lamperti type result Proposition \ref{prop0} holds. Taking into account the terminology, proposed in \cite{Paul20}, this r.f. has zero memory in both directions.

Now let us look at the case where the limit in the relation (\ref{lim1})  depends on the way how $(n_1, n_2)$ tends to infinity, i.e., the case where  we have the scaling transition. Let $n_1=n, n_2=n^\t$ and  let us consider the case $1/\al <\g_2\le \g_1<1$. Then $z_{\g_1,n_1}=n^{1-\g_1}, z_{\g_2,n_2}=n^{\t(1-\g_2)}$, and we define
$$
M_n=M_n(\g_1, \g_2):=\max\{z_{\g_1,n_1},z_{\g_2,n_2}\}=\max\{n^{1-\g_1} ,n^{\tau (1-\g_2)}\}=\begin{cases}
n^{1-\g_1}, \text{ if } \tau \leq \t_0,\\
n^{\tau (1-\g_2)}, \text{ if } \tau  > \t_0,
\end{cases}
$$
where $0<\t_0=(1-\g_1)/(1-\g_2)<1$. Finally, let us define
$$K_{\g_1}(\t)=\lim_{n\rightarrow\infty}\frac{z_{\g_1,n_1}}{M_n}=\begin{cases}
1, \text{ if } \tau \leq \t_0,\\
0, \text{ if } \tau > \t_0,
\end{cases} \quad
K_{\g_2}(\t)=\lim_{n\rightarrow\infty}\frac{z_{\g_2,n_2}}{M_n}=\begin{cases}
0, \text{ if } \tau < \t_0,\\
1, \text{ if } \tau \geq \t_0.
\end{cases}$$

\begin{prop}\label{prop3} Suppose that we have a sum (\ref{fieldd2}) of values of a linear r.f. with the filter (\ref{cij}). If $1/\al <\g_2\le \g_1<1 $ and $n_1=n, n_2=n^\t$, then, as $n \to \infty$,

\begin{equation}\label{lim2}
(n^{1+\t})^{-1/\al}M_n^{-1}S_{n,n^\t}(t_1,t_2)\fdd V(\t, t_1, t_2)
\end{equation}
where
\begin{equation}\label{lim2a}
V(\t, t_1, t_2)=\int_{-\infty}^{\infty}\int_{-\infty}^{\infty}\left( \ind{\{0\leq u_2\leq  t_2\}}K_{\g_1}(\t)H_{\g_1}(u_1,t_1)+
\ind{\{0\leq u_1\leq  t_1\}}K_{\g_2}(\t)H_{\g_2}(u_2,t_2)  \right)M(\dd u_1,\dd u_2).
\end{equation}
\end{prop}

From (\ref{lim2a}) we see, that linear r.f. $\{X_{k_1, k_2}\}$ in this case exhibits the scaling transition. Namely,  in the interval $0<\t<\t_0$ we have $K_{\g_1}(\t)=1 , \ K_{\g_2}(\t)=0$,  therefore  from (\ref{lim2a}) we get unbalanced scaling limit $V_-(t_1, t_2)$, independent of $\t$,  similarly in the interval $\t_0<\t<\infty$ we get unbalanced scaling limit $V_+(t_1, t_2)$. For $\t=\t_0$ we have  $K_{\g_1}(\t_0)= K_{\g_2}(\t_0)=1$ and we get balanced scaling limit $V(\t_0, t_1, t_2)$. If $1/\al<\g_1=\g_2=\g<1$ then we have $\t_0=1$ as a point of scaling transition and taking $n_1=n_2$ we have balanced scaling limit $V(1, t_1, t_2)$.

Also from this proposition one can see that there will be no scaling transition if only one  of  $\g_i, \ i=1, 2,$ does not satisfy the condition of Proposition \ref{prop2}, let us say that $\g_1> \max (1, 1/\al)$ and $1/\al<\g_2<1$. Then it is easy to see that  $z_{\g_1,n_1}=1$ and the normalization $z_{\g_2,n_2}$ for the second term in (\ref{fn1n2}) is growing to infinity, thus  it is prevailing. In this case the normalizing constant for $S_{n_1,n_2}(t_1, t_2)$ is $A_{n_1, n_2}=(n_1n_2)^{1/\al}z_{\g_2,n_2}$, and the limit process, independent of $\t$, will be obtained from (\ref{lim2a}) putting $K_{\g_1}(\t)\equiv 0 , \ K_{\g_2}(\t)\equiv 1$. In terms of directional memory one can say and the r.f. has zero memory in the horizontal direction and positive memory in the vertical direction.

\medskip

  Note that the linear r.f. from this example shows stronger effect than scaling transition as defined in Definition \ref{sctran}. Let  us consider the scale transition point $\t_0$  and balanced scaling limit process $V(\t_0, t_1, t_2)$. Let us take  $n_1=n, n_2=cn^{\t_0}$ with $0<c<\infty$, then we get $z_{\g_1,n_1}=n^{1-\g_1}, z_{\g_2,n_2}=n^{\t_0(1-\g_2)}c^{1-\g_2}$ and taking $M_n=n^{1-\g_1}$ we shall get
\begin{equation}\label{lim3}
(n^{-(1+\t_0)/\al+1-\g_1})c^{-1/\al}S_{n,cn^{\t_0}}(t_1, t_2)\fdd U(c, t_1, t_2)
\end{equation}
where
\begin{equation}\label{lim3a}
U(c, t_1, t_2)=\int_{-\infty}^{\infty}\int_{-\infty}^{\infty}\left( \ind{\{0\leq u_2\leq  t_2\}} H_{\g_1}(u_1,t_1)+
\ind{\{0\leq u_1\leq  t_1\}}c^{1-\g_2}H_{\g_2}(u_2,t_2)  \right)M(\dd u_1,\dd u_2).
\end{equation}
Taking  $M_n=n^{1-\g_1}\max (1, c^{1-\g_2})$ we can get as a limit the  process $V(\t_0, t_1, t_2)$, but now with both functions $K_{\g_i}$ depending on $c$.

This simple example gives us one more reason to reconsider Definition \ref{sctran}. Let us take in this example  $1/\al < \g_1<1$ and $\g_2=1$, then we get $z_{\g_1,n_1}=n_1^{1-\g_1}, z_{\g_2,n_2}=\ln n_2 $. If we assume, as earlier, $n_1=n, n_2=n^\t$, then $z_{\g_2,n_2}=\t \ln n$ and, for any $0<\t<\infty$, for sufficiently large $n$ we get $M_n(\g_1, \g_2):=\max\{z_{\g_1,n_1},z_{\g_2,n_2}\}=n^{1-\g_1}$. This would mean that there is no scale transition, but it is easy to see that such conclusion is due to our assumption about the relation between the growth of $n_2$ as  a function of $n$ (we  recall that $n_1=n$), namely, that $n_2=n^\t$. If we assume different relation, taking $n_2=\exp (n^\t)$, then  $z_{\g_2,n_2}=n^\t$ and we easily get the point of scale transition $\t_0=1-\g_1$. At the point of scale transition we have the same effect, which we described above and which can be named as the second order scale transition. It we assume $n_2=\exp (cn^{\t_0})$, then  it is easy to see that we shall get limit distribution dependent on the new parameter $c$. These considerations lead to the following generalization of Definition \ref{sctran}.
We consider a stationary r.f. $Y=\{Y_{k_1, k_2}, \ (k_1, k_2)\in \bz^2\}$ and sums defined in (\ref{sum}), only now we assume that $n_1=n, n_2=f(n, \t),$ where $\t\in (a, b)\subset \RR$ is some real parameter from an interval $(a, b)$, which can be finite or infinite, and $f: \ZZ_+ \times (a, b) \to \ZZ_+$. We suppose that for each fixed $\t$ function $f$ is monotonically growing to infinity, as $n\to \infty$. We denote
\begin{equation}\label{sumtau}
Z_{n, f,  \t}(t_1,t_2)=S_{n, f(n, \t)}(t_1,t_2), \quad  t_1\ge 0, \quad t_2\ge 0.
\end{equation}
We assume  that, for any $\t \in (a, b)$, there exists a nontrivial random field $V_{\t, f} (t_1,t_2)$ and a normalization $A_n (\t, f) \to \infty$ such that f.d.d. of $A_n^{-1}(\t, f) Z_{n, f, \t}(t_1,t_2)$ converges weakly to f.d.d. of $V_{\t, f}(t_1,t_2)$.

\begin{definition}\label{gensctran}
We say that  a random field $Y$ exhibits scaling transition if there exists a function $f(n, \t)$ and a point $\t_0\in (a, b)$ such that the limit process $V_{\t, f}$ is the same, let say $V_+$, for all  $\t \in (\t_0, \t_0+\d)$ and  another, not obtained by simple scaling, $V_-$, for   $\t \in (\t_0-\d, \t_0) $. The r.f. $V_{\t_0, f}$ is called well-balanced scaling limit of $Y$ at the point $\t_0$.

%We say that at the point $\t_0$ there is the second order scale transition if,  taking $n_2=f(cn, \t_0), c>0$ we get that the limit r.f. $V_{\t_0, f, c}$ depends on $c$.
\end{definition}
This definition not only extends the relation between $n_1$ and $n_2$ from power functions to more general class of functions, but also presupposes the possibility of more than one scale transition point, and such possibility will be realized in Example 3.

\subsection{Example 2}

We shall make very small change in the Example 1 in order to show  that  long-range dependence is not the main factor causing the scaling transition. In the filter (\ref{cij})  we redefine only $c_{0,0}:=a_1(0)+a_2(0)$, and $c_{0,0}$ is chosen in a such way that the following condition
\begin{equation}\label{sumZero}
\sum_{i_1=0}^{\infty}\sum_{i_2=0}^{\infty}c_{i_1,i_2}=0.
\end{equation}
is satisfied. Quantities $a_1(0), a_2(0)$ we choose in such a way, that $\sum_{i=0}^{\infty}a_j(i)=0, \ j=1, 2.$ Now the filter can be written as $c_{i_1,i_2}=a_1(i_1)\ind{i_1\geq 0}\ind{i_2=0}+a_2(i_2)\ind{i_2\geq 0}\ind{i_1=0}$, and, as in Example 1, we have formulae (\ref{chfuncd2})-(\ref{stgama})

Now we assume $a_j(i)=(1+i)^{-\g_j}, \ i\ge 1$,\ $\g_j>\max (1, 1/\al), \ j=1, 2,$ and, using the same notation (\ref{stgama}), we  investigate the  quantity (\ref{Jbnd2}).  We recall  that $a_0+b_0$ is negative and condition (\ref{sumZero}) holds.
 The integral in the expression of $J_{n_1, n_2}$ can be divided into three parts:
 \begin{equation*}
	\frac{J_{n_1, n_2}}{n_1n_2}= \int_{0}^{\infty}\int_{0}^{\infty}f_{n_1, n_2}^{(1)}(u_1,u_2)\dd u_1\dd u_2 + \int_{0}^{\infty}\left( \int_{-\infty}^{0}f_{n_1, n_2}^{(1)}(u_1,u_2)\dd u_1 \right)\dd u_2+ \int_{-\infty}^{0}\left( \int_{0}^{\infty}f_{n_1, n_2}^{(1)}(u_1,u_2)\dd u_1 \right)\dd u_2
\end{equation*}
\begin{equation*}	
	=J_{n_1, n_2}^{(1)}+J_{n_1, n_2}^{(2)}+J_{n_1, n_2}^{(3)}.
\end{equation*}
It is easy to see that
\begin{equation*}
	J_{n_1,n_2}^{(2)}=\int_{0}^{\infty}\left( \int_{-\infty}^{0}
	\abs{\sum_{l=1}^{d}x_l\left( \ind{\{0\leq n_2u_2\leq  \sv{n_2t_2^{(l)}}\}}
	U_{t_1^{(l)},\g_1}(\sv{n_1u_1},n_1)
		\right) }^\al
	\dd u_1 \right)
	 \dd u_2.
\end{equation*}
and
\begin{equation*}
	J_{n_1,n_2}^{(3)}=
	\int_{-\infty}^{0}
	\left( \int_{0}^{\infty}
	\abs{\sum_{l=1}^{d}x_l\left( \ind{\{0\leq n_1u_1\leq  \sv{n_1t_1^{(l)}}\}}U_{t_2^{(l)},\g_2}(\sv{n_2u_2},n_2)
		\right) }^\al
	\dd u_1 \right)
	 \dd u_2.
\end{equation*}

Since for sums, formed separately from sequences $\{a_1(k)\}$ or  $\{a_2(k)\}$ we have condition $\sum_{i=0}^{\infty}a_j(i)=0, \ j=1, 2$, therefore from Proposition \ref{prop1}
 we get
\begin{equation*}
	\frac{J_{n_1, n_2}^{(2)}}{n_1^{1-\g_1}}\rightarrow
	\int_{0}^{\infty}
	\left( \int_{-\infty}^{0}
	\abs{\sum_{l=1}^{m}x_l \ind{\{0\leq u_2\leq  t_2^{(l)}\}}
		H_{\g_1}(u_1,t_1^{(l)}) }^\al
	\dd u_1 \right)
	\dd u_2,
\end{equation*}
and
\begin{equation*}
\frac{J_{n_1, n_2}^{(3)}}{n_2^{1-\g_2}}\rightarrow
\int_{-\infty}^{0}\left( \int_{0}^{\infty}
\abs{\sum_{l=1}^{m}x_l \ind{\{0\leq u_1\leq  t_1^{(l)}\}}
	H_{\g_2}(u_2,t_2^{(l)}) }^\al
\dd u_1 \right)
\dd u_2.
\end{equation*}
It remains to consider the term $J_{n_1, n_2}^{(1)}$. Since $U_{t_1^{(l)},\g_1}(\sv{n_1u_1},n_1)=U_{t_1^{(l)},\g_1}(\sv{n_1u_1},n_1)\ind{\{0\leq \sv{n_1u_1}\leq  \sv{n_1t_1^{(l)}}\}}$ and similar equality can be written for $U_{t_2^{(l)},\g_2}(\sv{n_2u_2},n_2) $, we have
\begin{eqnarray}\label{jnm3}
J_{n_1, n_2}^{(1)} &=& \int_{\RR_+^2}\abs{\sum_{l=1}^{m}x_l\left( \ind{\{0\leq n_2u_2\leq  \sv{n_2t_2^{(l)}}\}}
	U_{t_1^{(l)},\g_1}(\sv{n_1u_1},n_1)+\ind{\{0\leq n_1u_1\leq  \sv{n_1t_1^{(l)}}\}}U_{t_2^{(l)},\g_2}(\sv{n_2u_2},n_2)  \right) }^\al\dd u_1\dd u_2 \\ \nonumber
 &=& \int_{\RR_+^2}
\abs{\sum_{l=1}^{m}x_l\ind{\{0\leq n_2u_2\leq  \sv{n_2t_2^{(l)}}\}}\ind{\{0\leq n_1u_1\leq  \sv{n_1t_1^{(l)}}\}}\left(
	U_{t_1^{(l)},\g_1}(\sv{n_1u_1},n_1)
	+	
	U_{t_2^{(l)},\g_2}(\sv{n_2u_2},n_2)  \right) }^\al\dd u_1\dd u_2.
\end{eqnarray}
Denoting
\begin{equation*}
U'_{t_j^{(l)},\g_j}(i,n_j):= \sum_{k=\sv{n_jt_j^{(l)}}-i+1}^{\infty}a_j(k),\ j=1, 2,
\end{equation*}
and taking into account (\ref{sumZero}), we have
\begin{equation*}
U_{t_1^{(l)},\g_1}(\sv{n_1u_1},n_1)
+	
U_{t_2^{(l)},\g_2}(\sv{n_2u_2},m)=-
\left( U'_{t_1^{(l)},\g_1}(\sv{nu},n)
+	
U'_{t_2^{(l)},\g_2}(\sv{n_2u_2},m) \right).
\end{equation*}
Substituting this equality into (\ref{jnm3}) we get
\begin{equation}\label{jnm4}
J_{n_1, n_2}^{(1)}=
\int_{0}^{\infty}\int_{0}^{\infty}
\abs{\sum_{l=1}^{m}x_l\ind{\{0\leq n_2u_2\leq  \sv{n_2t_2^{(l)}}\}}\ind{\{0\leq n_1u_1\leq  \sv{n_1t_1^{(l)}}\}}\left(
	U'_{t_1^{(l)},\g_1}(\sv{n_1u_1},n_1)
	+	
	U'_{t_2^{(l)},\g_2}(\sv{n_2u_2},n_2)  \right) }^\al\dd u_1\dd u_2.
\end{equation}
Now, as in Proposition \ref{prop1}, we can prove that
\begin{equation}\label{st1}
\frac{U'_{t_1^{(l)},\g_1}(\sv{n_1u_1},n_1)}{n_1^{1-\g_1}} \to \ind{\{0\leq u_1\leq  t_1^{(l)}\}}\frac{(t_1^{(l)}-u_1)^{1-\g_1}}{\g_1-1},
\end{equation}
\begin{equation}\label{st2}
\frac{U'_{t_2^{(l)},\g_2}(\sv{n_2u_2},n_2)}{n_2^{1-\g_2}} \to \ind{\{0\leq u_2\leq  t_2^{(l)}\}}\frac{(t_2^{(l)}-u_2)^{1-\g_2}}{\g_2-1},
\end{equation}
as $n_1, n_2\to \infty$.
Thus, we have the same situation as in Example 1 with sequences $z_{\g_1,n_1}=n_1^{1-\g_1}, z_{\g_2,n_2}=n_2^{1-\g_2}$, only now both sequences tend to zero, since $\g_i>1$. Taking $n_1=n, n_2=n^\t$ and
$$
M_n=M_n(\g_1, \g_2):=\max\{z_{\g_1,n_1},z_{\g_2,n_2}\}=\begin{cases}
n^{1-\g_1}, \text{ if } \tau \leq \t_0,\\
n^{\tau (1-\g_2)}, \text{ if } \tau  > \t_0,
\end{cases}
$$
where $0<\t_0=(\g_1-1)/(\g_2-1),$ we easily get that there is scale transition at the point $\t_0$.  It is possible to formulate the analog of Proposition \ref{prop3} with limit distributions for the appropriate normalized sum $S_{n,n^\t}(t_1,t_2)$ (to this aim we need to find majorizing functions as in (\ref{bound1}) and to use Lebesgue dominated convergence theorem), but we shall not do this, since the main task for this example was to show that scaling transition can be caused not by long-range dependence, which was present in the Example 1. In this example with the scale transition we have   condition (\ref{sumZero}), which indicates the negative dependence, this term was introduced in \cite{Lahiri} . Unfortunately, in both examples, in which we face the scale transition, it is not possible to use the notion of directional memory, which was introduced in \cite{Paul20}, since in both examples the normalizing constants are of the form
$$
A_{n_1, n_2}=(n_1n_2)^{1/\al}\max (n_1^{1-\g_1}, n_2^{1-\g_2}),
$$
and cannot be written in the form
$$
n_1^{1/\al+\d_1}n_2^{1/\al+\d_2}.
$$
But in the first example we can speak about general (not directional) positive memory, since the additional factor to $(n_1n_2)^{1/\al}$ is of the form $\max (n_1^{\d_1}, n_2^{\d_2})$ with both $\d_i$ positive, while in the second example these exponents $\d_i$ both are negative, therefore it is reasonable to speak about general negative memory.  It is possible to say that in the first example linear processes generated separately by filters $\{a_i\}$ and $\{b_i\}$ have positive memory, while in the second example, with conditions $\sum_{i=0}^\infty a_j(i)\ne 0, \ j=1, 2$,  both linear processes generated by these filters have negative memory.

\subsection{Example 3}

Now we shall take an example of a linear r.f. (\ref{fieldd2}) with the following filter
\begin{equation*}
c_{i,j}=
\begin{cases}
a_i, \text{ if } i\geq 1, j=0,\\
%b_j, \text{ if } j\geq 1, i=0,\\
c_i, \text{ if } i=j\geq 1,\\
0, \text{ otherwise},
\end{cases}
\end{equation*}
where $a_0=c_0=0, $ $a_i=(1+i)^{-\g_1}$, \ $c_i=(1+i)^{-\g_2}$, $1/\al<\g_2<\g_1<1$. It is convenient to write this filter as follows:
\begin{equation}\label{cij2}
c_{i,j}=a_{i}\ind{i\geq 1}\ind{j=0}+c_{j}\ind{i=j\geq 1},
\end{equation}
i.e., again we have coefficients of the filter on two lines, as in previous examples, but now  only one line is coordinate axis.

 As in Example 1, we assume  $n_1=n, n_2=n^\t$.  Now we get two points where scaling transition occurs. The first point $\t_0=(1-\g_1)/(1-\g_2)$ is the same as in Example 1, and we get additional point $\t_1=1$. Let us denote
\begin{equation}\label{An2}
A_n(\tau)=A_{n, n^\t}:=\begin{cases}
 n^{(1+\tau)/\al+1-\g_1} , \text{ if } \tau\in \left( 0, \t_0 \right),\\
 n^{(1+\tau)/\al+\tau(1-\g_2)} , \text{ if } \tau\in \left[\t_0  ,1 \right],\\
 n^{(1+\tau)/\al+1-\g_2} , \text{ if } \tau\in \left( 1,\infty \right),
\end{cases}
\end{equation}
and $a(u_1, u_2)=\max (0,-u_1,-u_2), \ b(u_1, u_2; t_1, t_2)=\min (t_1-u_1,t_2-u_2).$
Function $H_{\g_1}(u,t)$ was defined in Proposition \ref{prop1}, we recall here that for $\g<1$
\begin{equation*}
H_{\g}(u,t)=\frac{1}{1-\g}\left (((t-u)_+)^{1-\g}-((-u)_+)^{1-\g} \right ).
\end{equation*}

\begin{prop}\label{prop5} Suppose that we have a sum (\ref{fieldd2}) of values of a linear r.f. with the filter (\ref{cij2}). If $1/\al <\g_2\le \g_1<1 $ and $n_1=n, n_2=n^\t$, then, as $n \to \infty$,

\begin{equation}\label{lim4}
A_n^{-1}(\t)S_{n,n^\t}(t_1,t_2)\fdd V_1(\t, t_1, t_2)
\end{equation}
where
\begin{equation}\label{lim4a}
V_1(\t, t_1, t_2)=\begin{cases}
\int_{-\infty}^{\infty}\int_{-\infty}^{\infty}H_{\g_1}(u_1,t_1))\ind{\{0\leq u_2 \leq  t_2\}} M(\dd u_1 \dd u_2) , \text{ if } \tau\in \left( 0, \t_0 \right),\\
 \int_{-\infty}^{\infty}\int_{-\infty}^{\infty}\ind{\{0\leq u_1 \leq  t_1\}}H_{\g_2}(u_2,t_2)M(\dd u_1 \dd u_2)  , \text{ if } \tau\in \left(\t_0  ,1 \right),\\
 \int_{-\infty}^{\infty}\int_{-\infty}^{\infty} H_{\g_2}(u_1,t_1)\ind{\{0\leq u_2 \leq  t_2\}} M(\dd u_1 \dd u_2), \text{ if } \tau\in \left( 1,\infty \right).
\end{cases}
\end{equation}
and
\begin{equation}\label{lim4b}
V_1(\t_0, t_1, t_2)=\int_{-\infty}^{\infty}\int_{\infty}^{\infty}\left (\ind{\{0\leq u_2 \leq  t_2\}}H_{\g_1}(u_1,t_1)+ \ind{\{0\leq u_1 \leq  t_1\}}H_{\g_2}(u_2,t_2)\right ) M(\dd u_1 \dd u_2),
\end{equation}

\begin{equation}\label{lim4c}
V_1(1, t_1, t_2)=\int_{-\infty}^{\infty}\int_{\infty}^{\infty}\ind{a(u_1, u_2) <b(u_1, u_2; t_1, t_2)}\frac{(b(u_1, u_2; t_1, t_2))^{1-\g_2}-(a(u_1, u_2))^{1-\g_2}}{1-\g_2} M(\dd u_1 \dd u_2).
\end{equation}
\end{prop}

The sum $Z_{\bt^{(l)},\g_2}(i,j,n_1,n_2)$   is more complicated comparing with $U_{t_j^{(l)},\g_j}(i,n_j)$ (see (\ref{stgama}), and it alone gives us one point of transition.
Namely, we can consider r.f. (\ref{fieldd2}) with the following simple filter
\begin{equation}\label{cijlygus}
c_{i,j}=c_{i}, \ \ {\rm if} \  i=j\ge 1, \ \ {\rm and} \ \ c_{i,j}=0, \ \ {\rm elsewhere}.
\end{equation}
Analysis of this random field gives us the scaling transition point $\t=1$.

{\it Sketch of the proof of Proposition \ref{prop5}}.
Taking into account formulae (\ref{chfuncd2})-(\ref{fn1n2}) and considering  the ch.f. of the vector $\left( S_\bn(\bt^{(1)}),\dots, S_\bn(\bt^{(m)}) \right)$, we must investigate the asymptotic of $f_{n_1,n_2}^{(3)}(u_1,u_2)$. Similarly to Example 1, taking into account (\ref{cij2}), we can get the following expression of this quantity
\begin{equation} \label{CijnmIsraiska}
f_{n_1,n_2}^{(3)}(u_1,u_2) =\abs{ \sum_{l=1}^{m}x_l\left(\ind{\{0\leq n_2u_2\leq  \sv{n_2t_2^{(l)}}\}}
U_{t_1^{(l)},\g_1}(\sv{n_1u_1},n_1)+Z_{\bt^{(l)},\g_2}(\sv{n_1u_1}, \sv{n_2u_2},n_1,n_2)
 \right)}^\al,
\end{equation}
where $U_{t_j^{(l)},\g_j}(i,n)$ is defined in (\ref{stgama}) and
\begin{equation}\label{ztsl}
Z_{\bt^{(l)},\g_2}(i,j,n_1,n_2):=\sum_{k=\max\left( 0,-i,-j \right)}^{\min\left( \sv{n_1t_1^{(l)}}-i,\sv{n_2t_2^{(l)}}-j \right)}c_{k}.
\end{equation}
Comparing the expression (\ref{CijnmIsraiska}) with the corresponding expression of $f_{n_1,n_2}(u_1,u_2)$ in Example 1 we see that the first term in (\ref{CijnmIsraiska}) is the same (corresponding to the filter on the horizontal axis), and for this term we can apply Proposition \ref{prop1}. This will give us the following relation (independent of $\t$, since $U_{t_1^{(l)},\g_1}(\sv{n_1u_1},n_1)$ depends only on $n_1$):
\begin{equation*}
\frac{U_{t_1^{(l)},\g_1}(\sv{nu_1},n_1)}{z_{\g_1,n_1}}=\rightarrow
\ind{(-\infty,t_l)}(u_1)\int_{(-u_1)_+}^{t_l-u_1}v^{-\g_1} \dd v =H_{\g_1}(u_1,t_l),
\end{equation*}
with $z_{\g_1,n_1}=n_1^{1-\g_1}.$
The second term (corresponding to the filter on the diagonal) in (\ref{CijnmIsraiska}) is different, and Proposition \ref{prop1} directly cannot be applied. Therefore, we investigate  the quantity
\begin{equation}\label{jnm2}
J_{n_1,n_2}=n_1n_2\int_{-\infty}^{\infty}\int_{-\infty}^{\infty}f_{n_1,n_2}^{(3)}(u_1,u_2){\dd u}_1{\dd u}_2=J^{(1)}_{n_1,n_2}+J^{(2)}_{n_1,n_2}+J^{(3)}_{n_1,n_2}+J^{(4)}_{n_1,n_2}
\end{equation}
with $n_1=n, n_2=n^\t$ (but sometimes we shall leave $n_1, n_2$ instead of $n, n^\t$), dividing the integral over $\RR^2$ into four integrals over regions $\{u\ge 0, v\ge 0 \}, \{u< 0, v\ge 0 \}, \{u< 0, v< 0 \}, \{u\ge 0, v< 0 \}$. In investigation of these integrals we use without special mentioning the following steps: we prove the point-wise convergence of functions in the expression of integrals and then use the Lebesgue dominated convergence theorem.

(i) We start with the integral
$$
J^{(1)}_{n,n^\t}=n^{1+\t}\int_{0}^{\infty}\int_{0}^{\infty}f_{n,n^\t}^{(3)}(u_1,u_2){\dd u}_1{\dd u}_2
$$
and consider the case $0<\t<1$.Changing the sum in (\ref{ztsl}) by integral after some transformations we can get
$$
Z_{\bt^{(l)},\g_2}(\sv{n_1u_1}, \sv{n_2u_2}, n_1, n_2)= n_2^{1-\g_2}\ind{\{\sv{n_1u_1}\leq \sv{n_1t_1^{(l)}}\}}\ind{\{\sv{n_2u_2}\leq \sv{n_2t_2^{(l)}}\}}
\int_{0}^{\infty}\kappa_{n_1,n_2}(y)\dd y,
$$
where
$$
\kappa_{n_1,n_2}(y)=\ind{\left( 0, \min\left( (\sv{n_1t_1^{(l)}}-\sv{n_1u_1}+1)/{n_2},(\sv{n_2t_2^{(l)}}-\sv{n_2u_2}+1)/{n_2} \right) \right)}(y)
n_2^{\g_2}c_{\sv{n_2y}}.
$$
Since for a fixed $y>0$ we have
$\kappa_{n,m}(y)\rightarrow \ind{\left( 0, s_l-u_2 \right) }(y) y^{-\g_2}$,
and this function can be bounded by integrable function $\abs{\kappa_{n,m}(y)}\leq
\ind{\left( 0, s_l-u_2+2 \right) }(y) y^{-\g_2},$ we get
$$
\int_{0}^{\infty}\kappa_{n_1,n_2}(y)\dd y\rightarrow
\int_{0}^{t_2^{(l)}-u_2} y^{-\g_2} \dd y=\frac{(t_2^{(l)}-u_2)^{1-\g_2}}{1-\g_2}.
$$
Then, for $u_1\ge0, u_2\ge 0, $ we can get
$$
\frac{Z_{\bt^{(l)},\g_2}(\sv{n_1u_1}, \sv{n_2u_2}, n_1, n_2)}{n_2^{1-\g_2}}\rightarrow
\ind{\{u_1 <t_1^{(l)}\}}
\ind{\{u_2 <t_2^{(l)}\}}
\frac{(t_2^{(l)}-u_2)^{1-\g_2}}{1-\g_2}=\ind{\{u_1 <t_1^{(l)}\}} H_{\g_2}(u_2,t_2^{(l)}),
$$
since, for $u_2\ge 0,$ we have $(-u_2)_+=0$. Now we return to the function from (\ref{CijnmIsraiska})
$$
f_{n_1,n_2}^{(3)}(u_1,u_2)  = \Big | \sum_{l=1}^{m}x_l\Big (\ind{\{0\leq n_2u_2\leq  \sv{n_2t_2^{(l)}}\}}
n_1^{1-\g_1}\frac{U_{t_1^{(l)},\g_1}(\sv{n_1u_1},n_1)}{n_1^{1-\g_1}}
$$
$$
+  n_2^{1-\g_2}\frac{Z_{\bt^{(l)},\g_2}(\sv{n_1t_1^{(l)}}, \sv{n_2t_2^{(l)}},n_1,n_2)}{n_2^{1-\g_2}}
 \Big )\Big |^\al,
$$
Thus we have the same situation as in Example 1. Denoting $M_n:=\max(n^{1-\g_1},m^{1-\g_2})=\max(n^{1-\g_1},n^{\tau(1-\g_2)})$ we get the same functions $K_i(\t)$:
$$
\lim_{n\to \infty}\frac{n^{1-\g_1}}{M_n}=
K_1(\tau):=\begin{cases}
0,\text{ if }1>\tau >\t_0,\\
1,\text{ if }0<\tau \leq \t_0,
\end{cases}
$$
$$
\lim_{n\to \infty}\frac{m^{1-\g_2}}{M_n}=
K_2(\tau):=\begin{cases}
1,\text{ if }\tau \geq \t_0,\\
0,\text{ if }\tau < \t_0.
\end{cases}
$$
Then we easily get
$$
\frac{f_{n_1,n_2}^{(3)}(u_1,u_2)}{M_n^\al}\rightarrow\Big | \sum_{l=1}^{m}x_l\Big (\ind{\{0\leq u_2 \leq  t_2^{(l)}\}}
K_1(\tau) H_{\g_1}(u_1,t_1^{(l)}) +\ind{\{u_1 <t_1^{(l)}\}}K_2(\tau) H_{\g_2}(u_2,t_2^{(l)})\Big )\Big |^\al
$$
and
\begin{equation}\label{JNM1}
\frac{J^{(1)}_{n,n^\t}}{n^{1+\t}M_n^\al}\rightarrow \int_{0}^{\infty}\int_{0}^{\infty}\Big | \sum_{l=1}^{m}x_l\Big (\ind{\{0\leq u_2 \leq  t_2^{(l)}\}}
K_1(\tau) H_{\g_1}(u_1,t_1^{(l)})+\ind{\{u_1 <t_1^{(l)}\}}K_2(\tau) H_{\g_2}(u_2,t_2^{(l)})\Big )\Big |^\al {\dd u}_1{\dd u}_2.
\end{equation}
In the case $\t>1$ (this means $n_2/n_1\to \infty$) we similarly can get (we recall that $u_1\ge0, u_2\ge 0 $)
$$
\frac{Z_{\bt^{(l)},\g_2}(\sv{n_1u_1}, \sv{n_2u_2}, n_1, n_2)}{n_1^{1-\g_2}}\rightarrow
\ind{\{u_1 <t_1^{(l)}\}}
\ind{\{u_2 <t_2^{(l)}\}}
\frac{(t_1^{(l)}-u_1)^{1-\g_2}}{1-\g_2}=\ind{\{u_2 <t_2^{(l)}\}} H_{\g_2}(u_1,t_1^{(l)}).
$$
Since for $n_1=n, n_2=n^\t$ norming constants for the two terms in (\ref{CijnmIsraiska}) are $n^{1-\g_1}$ and $n^{1-\g_2}$, respectively, and $n^{1-\g_1}<n^{1-\g_2}$, we get
\begin{equation}\label{JNM1A}
\frac{J^{(1)}_{n,n^\t}}{n^{1+\t+(1-\g_2)\al}}\rightarrow \int_{0}^{\infty}\int_{0}^{\infty}\Big | \sum_{l=1}^{m}x_l\Big (\ind{\{0\leq u_2 \leq  t_2^{(l)}\}}
 H_{\g_2}(u_1,t_1^{(l)}) \Big )\Big |^\al {\dd u}_1{\dd u}_2.
\end{equation}

(ii) Now we investigate
$$
J^{(2)}_{n,n^\t}=n^{1+\t}\int_{-\infty}^{0}\int_{0}^{\infty}f_{n,n^\t}^{(3)}(u_1,u_2){\dd u}_2{\dd u}_1
$$
and consider the case $0<\t<1$.Changing the sum in (\ref{ztsl}) by integral after some transformations we can get
$$
Z_{\bt^{(l)},\g_2}(\sv{n_1u_1}, \sv{n_2u_2}, n_1, n_2)\leq
\ind{\brc{-\sv{n_1u_1}\leq \sv{n_2t_2^{(l)}}-\sv{n_2u_2} }}
\ind{\brc{\sv{n_2u_2}\leq \sv{n_2t_2^{(l)}} }}
\int_{0}^{ {n_2t_2^{(l)}}-{n_2u_2} +2n_2}
y^{-\g_2} \dd y
$$
$$
=\ind{\brc{-\sv{n_1u_1}\leq \sv{n_2t_2^{(l)}}-\sv{n_2u_2} }}
\ind{\brc{\sv{n_2t_2}\leq \sv{n_2t_2^{(l)}} }}
\frac{\skl{ {n_2t_2^{(l)}}-{n_2u_2} +2n_2}^{1-\g_2}}{1-\g_2},
$$
hence, using simple inequalities between indicator functions we get
$$
\frac{Z_{t_1^{(l)},t_2^{(l)},\g_2}(\sv{n_1u_1},\sv{n_2t_2}, n_1, n_2)}{n_2^{1-\g_2}}\leq
\ind{\brc{u_1\geq \frac{n_2}{n_1}(t_2-t_2^{(l)}-1) }}
\ind{\brc{t_2\leq u_2^{(l)}+1 }}
\frac{\skl{ {t_2^{(l)}}-{u_2} +2}^{1-\g_2}
}{1-\g_2}.
$$
Since we consider the set $\{(u_1,u_2): u_1<0,\ u_2\geq 0\}$, we have $\ind{\brc{u_1\geq \frac{n_2}{n_1}(u_2-t_2^{(l)}-1) }}\rightarrow \ind{\brc{u_1\geq 0 }}= 0,$ therefore,
$$
\frac{Z_{t_1^{(l)},t_2^{(l)},\g_2}(\sv{n_1u_1},\sv{n_2t_2}, n_1, n_2)}{n_2^{1-\g_2}}\rightarrow 0.
$$
Using the same sequence $M_n$ and the function $K_1,$ we get
$$
\frac{f_{n,n^\t}^{(3)}(u_1,u_2)}{M_n^\al}\rightarrow \Big | \sum_{l=1}^{m}x_l \ind{\{0\leq u_2 \leq  t_2^{(l)}\}}
K_1(\tau) H_{\g_1}(u_1,t_1^{(l)})\Big |^\al,
$$
and, finally,
\begin{equation}\label{JNM2}
\frac{J^{(2)}_{n,n^\t}}{n^{1+\t}M_n^\al}\rightarrow \int_{-\infty}^{0}\int_{0}^{\infty}\Big | \sum_{l=1}^{m}x_l \ind{\{0\leq u_2 \leq  t_2^{(l)}\}}
K_1(\tau) H_{\g_1}(u_1,t_1^{(l)})\Big |^\al {\dd u}_2{\dd u}_1.
\end{equation}

In the case $\t>1$ we  get
$$
Z_{\bt^{(l)},\g_2}(\sv{n_1u_1}, \sv{n_2u_2}, n_1, n_2)=n_1^{1-\g_2}\ind{\{\sv{n_2u_2}\leq \sv{n_2t_2^{(l)}}\}}\int_{0}^{\infty}\kappa_{n_1,n_2}(y)\dd y,
$$
where
$$
\kappa_{n_1,n_2}(y)=\ind{\left( -\sv{n_1u_1}/n_1, \min\left( (\sv{n_1t_1^{(l)}}-\sv{n_1u_1}+1)/{n_1}, (\sv{n_2t_2^{(l)}}-\sv{n_2u_2}+1)/{n_1} \right) \right)}(y)
n_1^{\g_2}c_{\sv{n_2y}}.
$$
Using the relation $\kappa_{n_1,n_2}(y) \rightarrow \ind{\left( -u_1, t_1^{(l)}-u_1 \right) }(y) y^{-\g_2},$ in the same way as in (i) we obtain
$$
\frac{Z_{\bt^{(l)},\g_2}(\sv{n_1u_1}, \sv{n_2u_2}, n_1, n_2)}{n_1^{1-\g_2}}\rightarrow \ind{\{u_2 <t_2^{(l)}\}}H_{\g_2}(u_1,t_1^{(l)}).
$$
Again, due to $n^{1-\g_1}<n^{1-\g_2}$, the second term in (\ref{CijnmIsraiska}) is prevailing,  therefore we get
\begin{equation}\label{JNM2A}
\frac{J^{(2)}_{n,n^\t}}{n^{1+\t+(1-\g_2)\al}}\rightarrow \int_{-\infty}^{0}\int_{0}^{\infty}\Big | \sum_{l=1}^{m}x_l \ind{\{0\leq u_2 \leq  t_2^{(l)}\}}
 H_{\g_2}(u_1,t_1^{(l)})\Big |^\al {\dd u}_2{\dd u}_1.
\end{equation}

(iii) We investigate the third integral
$$
J^{(3)}_{n,n^\t}=n^{1+\t}\int_{-\infty}^{0}\int_{-\infty}^0 f_{n,n^\t}^{(3)}(u_1,u_2){\dd u}_2{\dd u}_1
$$
Since our goal is to prove that after the appropriate normalization this integral tends to zero, we use the rough estimate
$$
J^{(3)}_{n,n^\t}\leq n^{1+\t}d^\al\sum_{l=1}^{d}\abs{x_l}^\al\int_{-\infty}^{0}\int_{-\infty}^0\abs{Z_{\bt^{(l)},\g_2}(\sv{n_1u_1}, \sv{n_2u_2}, n_1, n_2)}^\al {\dd u}_2{\dd u}_1.
$$
Let us assume that $\t<1,$ i.e., $n_2/n_1 \to 0$. By change of variables we have
$$
\int_{-\infty}^{0}\int_{-\infty}^0\abs{Z_{\bt^{(l)},\g_2}(\sv{n_1u_1}, \sv{n_2u_2}, n_1, n_2)}^\al {\dd u}_2{\dd u}_1
$$
$$
=\int_{-\infty}^{0}\int_{-\infty}^{-\sv{n_2u_2}/n_1}\abs{Z_{\bt^{(l)},\g_2}(\sv{n_1u_1}+\sv{n_2u_2}, \sv{n_2u_2}, n_1, n_2)}^\al {\dd u}_1{\dd u}_2.
$$
Once more, changing the sum into integral, we get
$$
Z_{\bt^{(l)},\g_2}(\sv{n_1u_1}+\sv{n_2u_2}, \sv{n_2u_2}, n_1, n_2)
$$
$$
\leq \ind{\brc{\sv{n_1u_1}\leq \sv{n_1t_1^{(l)}}}}
\ind{\brc{\sv{n_1u_1}\geq -\sv{n_2t_2^{(l)}}}}
n_2^{1-\g_2}
\frac{\left( t_2^{(l)}-u_2+1 \right)^{1-\g_2}-\left( -u_2 \right)^{1-\g_2}}{1-\g_2}.
$$
This gives us the following estimate
$$
\frac{J^{(3)}_{n_1,n_2}}{n_1n_2^{1+\al(1-\g_2)}}\leq \int_{-\infty}^{0}|H_{\g_2}(u_2,t_2^{(l)}+1)|\int_{-\infty}^{\infty} \ind{(-\infty,-\sv{n_2u_2}/n_1)}(u_1)\ind{\brc{\sv{n_1u_1}\leq \sv{n_1t_1^{(l)}}}}
\ind{\brc{\sv{n_1u_1}\geq -\sv{n_2t_2^{(l)}}}}{\dd u}_1{\dd u}_2.
$$
Considering these three indicator functions we get
$$
\ind{(-\infty,-\sv{n_2u_2}/n_1)}(u_1)\ind{\brc{\sv{n_1u_1}\leq \sv{n_1t_1^{(l)}}}}
\ind{\brc{\sv{n_1u_1}\geq -\sv{n_2t_2^{(l)}}}}\leq \ind{\left[-\frac{n_2}{n_1}t_2^{(l)} , \min(t_1^{(l)}+1,-\sv{n_2u_2}/n_1 )\right]}(u_1),
$$
therefore,
$$
\int_{-\infty}^{\infty} \ind{(-\infty,-\sv{n_2u_2}/n_1)}(u_1)\ind{\brc{\sv{n_1u_1}\leq \sv{n_1t_1^{(l)}}}}
\ind{\brc{\sv{n_1u_1}\geq -\sv{n_2t_2^{(l)}}}}{\dd u}_1\leq \min(t_1^{(l)}+1,-\frac{\sv{n_2}}{n_1u_2} )+\frac{n_2}{n_1}t_2^{(l)}.
$$
For a fixed $u_2$ the right-hand side of the last inequality tends to zero, therefore, we get
\begin{equation}\label{JNM3}
\frac{J^{(3)}_{n_1,n_2}}{n_1n_2^{1+(1-\g_2)\al}}\rightarrow 0.
\end{equation}

In the case $\t>1,$ i.e., $n_1/n_2 \to 0$, we get
\begin{equation}\label{JNM3A}
\frac{J^{(3)}_{n_1,n_2}}{n_1^{1+(1-\g_2)\al}n_2}\rightarrow 0.
\end{equation}

(iv) It remains to investigate the last integral
$$
J^{(4)}_{n_1,n_2}=n_1n_2\int_{0}^{\infty}\int_{-\infty}^0 f_{n_1,n_2}^{(3)}(u_1,u_2){\dd u}_2{\dd u}_1
$$
$$
=n_1n_2\int_{0}^{\infty}\int_{-\infty}^0 \abs{Z_{\bt^{(l)},\g_2}(\sv{n_1u_1}, \sv{n_2u_2}, n_1, n_2)}^\al {\dd u}_2{\dd u}_1.
$$
Let us note that
$$
Z_{t_1^{(l)}, t_2^{(l)}, \g_2}(\sv{n_1u_1}, \sv{n_2u_2}, n_1, n_2)=Z_{t_2^{(l)}, t_1^{(l)}, \g_2}(\sv{n_2u_2}, \sv{n_1u_1},  n_2, n_1),
$$
therefore the investigation of this quantity in the case $\t<1$ will be the same as in  the case (ii) and $\t>1$. We shall get
\begin{equation}\label{JNM4}
\frac{J^{(4)}_{n_1,n_2}}{n_1n_2^{1+(1-\g_2)\al}}\rightarrow \int_{-\infty}^{0}\int_{0}^{\infty}\Big | \sum_{l=1}^{m}x_l \ind{\{0\leq u_2 \leq  t_2^{(l)}\}}
 H_{\g_2}(u_1,t_1^{(l)})\Big |^\al {\dd u}_2{\dd u}_1.
\end{equation}
In the case $\t>1$ investigation is similar to (ii) and $\t<1$, and we get
\begin{equation}\label{JNM4A}
\frac{J^{(4)}_{n_1,n_2}}{n_1^{1+(1-\g_2)\al}n_2}\rightarrow 0.
\end{equation}
Collecting formulae (\ref{JNM1}) - (\ref{JNM4A}) we get (\ref{lim4a}). Due to the fact that $K_1(\t_0)=K_2(\t_0)=1$ we get (\ref{lim4b}).

It is more difficult to verify (\ref{lim4c}), to this aim we must simply go through all the proof of (\ref{lim4a}), assuming that $n_1=n_2=n$ and $A_n=n^{(2/\al)+1-\g_2}$. Let us take the case (i), the integral $J^{(1)}_{n,n}$. In all four integrals the first term in the expression of the function $f_{n_1,n_2}^{(3)}(u_1,u_2)$ from (\ref{CijnmIsraiska}) is independent of $n_2$, therefore, we need to consider only the second term, namely, $Z_{t_1^{(l)}, t_2^{(l)}, \g_2}$. In the case $n_1=n_2=n$ it is easy to see that $\kappa_{n,n}(y)\rightarrow \ind{\left( 0, b(u_1, u_2; t_1^{(l)}, t_2^{(l)}) \right) }(y) y^{-\g_2}$, therefore, for $u_1\ge 0, u_2\ge 0,$  we get
$$
\frac{Z_{\bt^{(l)},\g_2}(\sv{nu_1}, \sv{nu_2}, n, n)}{n^{1-\g_2}}\rightarrow
\ind{\{u_1 <t_1^{(l)}\}}
\ind{\{u_2 <t_2^{(l)}\}}
\frac{b(u_1, u_2; t_1^{(l)}, t_2^{(l)})^{1-\g_2}}{1-\g_2}.
$$
Since $\g_2<\g_1$ we see that the first term after normalization tends to zero, also we have $a(u_1, u_2)=0,$ therefore we see that we got the same expression as in (\ref{lim4c}), in the case  $u_1\ge 0, u_2\ge 0.$ Similarly can be treated the rest three integrals in (\ref{jnm2}).
\halmos

\section{The case $d=3$}

In the above mentioned papers \cite{Puplinskaite2},  \cite{Puplinskaite1}, and \cite{Pilipauskaite} only  r.f.s, defined on $\ZZ^2$ were considered. It was mentioned that generalization of the results of these papers to the case $\ZZ^d$ with $d\ge 3$ is difficult task, and the first step in this direction was the paper \cite{Surg}, where the scaling transition of linear random fields on $\ZZ^3$ was considered. But at first we must discuss what we understand by words scaling transition on $\ZZ^3$, or more generally, on $\ZZ^d, \ d\ge 3$. In section 2 (the case $d=2$) considering sums (\ref{sum}) we have two possibilities. The first one is the convergence (in the sense of f.d.d.) of  appropriately normed sums   to some limit process as $(n_1, n_2)\to \infty$, and the limit process is independent on the way how the indices $n_1,n_2$ grow. The second one is the situation when the limit for sums (\ref{sum}) depends on the way how $(n_1, n_2)$ grow. In this case there  was quite natural way to define this dependence, using   relation $n_2=f(n_1)$ between $n_1$ and $n_2$, see Definitions 2 and 6. But the straightforward generalization of Definition 2, considering sums (\ref{field}) in the case $d=3$ and assuming $n_1=n^{q_1}, n_2=n^{q_2}, n_3=n^{q_3}$ (such case is considered in \cite{Surg}), is too narrow, since it presents only one possible way to define the path in $\ZZ_+^3$. Probably for this reason  in \cite{Surg} there is no strict definition of the scaling transition, and the author  in \cite{Surg} wrote "...we do not attempt to provide a formal definition of scaling transition for RFs in dimensions $d\ge 3$ since further studies are needed to fully understand it". Our examples of this section show that in dimension 3 we have much more possibilities (comparing with the case $d=2$) to define paths of $(n_1, n_2, n_3)$ growing to infinity, moreover, with growing dimension the complexity grows very rapidly. Therefore, we propose to define the scaling transition, independently of dimension $d$, as the case where the limits for sums (\ref{field}) depend on the way how $\bn \to \infty$, i.e., there exist at least two paths in $\ZZ_+^d$ such that limit r.f. for these paths are different and cannot be obtained one from another by simple scaling. With such definition  we should have a simple dichotomy in limit theorems for sums of values of random fields (\ref{field}): or the limit theorem of Lamperti type, as formulated in Corollary 1 in \cite{DavPau1} (see Proposition 1 in the case $d=2$), holds, either there is the scaling transition.

At first  we shall show  that it is  easy to generalize  examples of r.f. on $\ZZ^2$ with simple structure of filters considered in the previous subsections, to higher dimensions. We consider the case $d=3$, although generalization to higher dimensions in some examples does not present any principal difficulties. For the notation of three-dimensional multi-indices we  use bold letters, for example $\bi =(i_1, i_2, i_3), \bn =(n_1, n_2, n_3)$. We consider linear r.f.
\begin{equation}\label{field3d}
X_{\bk}=\sum_{\bi \in \ZZ^3_+} c_{\bi}\xi_{\bk-\bi}, \ \bk \in \bz^3,
\end{equation}
where $\ZZ^3_+=\{\bi \in \ZZ^3: \ \bi\ge {\B0} \}$, and investigate the asymptotic behavior of the process
\begin{equation}\label{sum13d}
 	S_{\bn}(\bt)=\sum_{\bk=\B0}^{\sv{\bn \bt}}X_{\bk}.
 \end{equation}

\subsection{Example 4}

 We take three sequences of positive numbers $a_j(i)=(1+i)^{-\g_j},  \ i\ge 1, \g_j>1/\al, \  a_j(0)=0, \ j=1, 2,3,$ and the following filter
\begin{equation}\label{ci3d}
c_{\bi}=
\begin{cases}
a_{1}(i_1), \text{ if } i_1\geq 1, i_2=i_3=0,\\
a_{2}(i_2), \text{ if } i_2\geq 1, i_1=i_3=0,\\
a_{3}(i_3), \text{ if } i_3\geq 1, i_2=i_1=0,\\
0, \text{ elsewhere}.
\end{cases}
\end{equation}
 Since this filter is exactly the same as considered in the section "Preliminaries" and in Example 1 (in the case $d=2$), we can start with formula (\ref{fn2}) and to write  formulae, similar to (\ref{Jbnd2}) and (\ref{fn1n2}):

\begin{equation}\label{Jbn}
J_{\bn}=n_1n_2n_3
\int_{-\infty}^{\infty}\int_{-\infty}^{\infty}\int_{-\infty}^{\infty}f_{\bn}^{(4)}(\bu)\dd u_1\dd u_2\dd u_3,
\end{equation}
where
\begin{equation}\label{fbn}
f_{\bn}^{(4)}(\bu)=\abs{\sum_{l=1}^{d}x_l\left( \sum_{p=1}^3 B_p(\bn, \bu, l) U_{t^{(l)}_p,\g_p}(\sv{n_pu_p},n_p)  \right) }^\al,
\end{equation}
\begin{equation}\label{Bpbn}
B_p(\bn, \bu, l)=\ind{\{0\leq \sv{n_ru_r} \leq \sv{n_rt^{(l)}_r}, \ 0\leq \sv{n_su_s} \leq \sv{n_st^{(l)}_s}, r\ne s, r, s \in Q_p\}}, \quad Q_p=\{1, 2, 3\}\setminus p,
\end{equation}
and
\begin{equation}\label{Sngamma}
U_{t^{(l)}_p,\g_p}(i_p,n_p)=\sum_{k=\left( -i_p \right)\vee 0 }^{\sv{n_pt^{(l)}_p}-i_p}a_{k}^{(p)}.
\end{equation}
For quantities $U_{t^{(l)}_p,\g_p}(i_p,n_p)$ we apply Proposition \ref{prop1} with normalization by quantities $z_{\g_p,n_p}^{(p)}, \ p=1, 2, 3,$ and the further analysis depends on our assumptions on exponents $\g_p$ and relation between coordinates of $\bn$. For example, assuming that all $\g_p>1$ we shall get that there is no scaling transition and we have the analogue of Proposition \ref{prop2}. Let us consider the case which will give us the scaling transition. We assume $1/\al<\g_1<\g_2<\g_3<1$ and  we take $n_1=n$ and , for some positive real numbers $\t$ and $\s$, we set $n_2=n^\t, n_3=n^\s$. Then the normalization quantities $z_{\g_p,n_p}^{(p)}, \ p=1, 2, 3,$ become functions of $n$:
$$
z_{\g_1,n}^{(1)}=n^{1-\g_1}, z_{\g_2,n}^{(2)}=n^{\t(1-\g_2)}, \ z_{\g_3,n}^{(3)}=n^{\s(1-\g_3)}.
$$
Applying Proposition \ref{prop1} we get, as $n\to \infty$,
\begin{equation}\label{lim5}
\frac{U_{t^{(l)}_p,\g_p}(\sv{n_pu_p},n_p)}{z_{\g_p,n}^{(p)}} \to H_{\g_p}(u_p,t^{(l)}_p).
\end{equation}
We define
$$
M_{n}(\t, \s):=\max\{z_{\g_p,n}^{(p)}, \ p=1, 2, 3 \}=\max\{n^{1-\g_1}, n^{\tau (1-\g_2)}, n^{\s(1-\g_3)}\}
$$
and
\begin{equation}\label{Kp}
K_p(\t,\s)=\lim_{n\rightarrow\infty}\frac{z_{\g_p,n}^{(p)}}{M_{n}(\t, \s)}.
\end{equation}
Let us denote
 $$
 \t_0=\frac{ 1-\g_1}{ 1-\g_2}, \ \s_0=\frac{ 1-\g_1}{ 1-\g_3}, \ a=\frac{ 1-\g_2}{ 1-\g_3}. %\ b=\frac{ 1-\g_1}{ 1-\g_3},
 $$

 \begin{figure}[!ht]
 	\begin{subfigure}[b]{0.33\linewidth}
 		\centering
 		% Created by tikzDevice version 0.12 on 2019-03-19 10:39:35
% !TEX encoding = UTF-8 Unicode
\begin{tikzpicture}[x=1pt,y=1pt]
\definecolor{fillColor}{RGB}{255,255,255}
\path[use as bounding box,fill=fillColor,fill opacity=0.00] (0,0) rectangle (144.54,144.54);
\begin{scope}
\path[clip] (  0.00,  0.00) rectangle (144.54,144.54);
\definecolor{drawColor}{RGB}{0,0,0}

\path[draw=drawColor,line width= 0.8pt,line join=round,line cap=round] ( 16.51, 16.51) -- (128.03, 16.51);

\path[draw=drawColor,line width= 0.8pt,line join=round,line cap=round] (124.64, 15.27) --
	(128.03, 16.51) --
	(124.64, 17.74);

\path[draw=drawColor,line width= 0.8pt,line join=round,line cap=round] ( 16.51, 16.51) -- ( 16.51,128.03);

\path[draw=drawColor,line width= 0.8pt,line join=round,line cap=round] ( 17.74,124.64) --
	( 16.51,128.03) --
	( 15.27,124.64);

\node[text=drawColor,anchor=base,inner sep=0pt, outer sep=0pt, scale=  1.00] at (131.38,  9.55) {$\tau$};

\node[text=drawColor,anchor=base,inner sep=0pt, outer sep=0pt, scale=  1.00] at ( 13.16,128.88) {$\sigma$};

\node[text=drawColor,anchor=base,inner sep=0pt, outer sep=0pt, scale=  1.00] at ( 53.68,  5.86) {1};

\node[text=drawColor,anchor=base,inner sep=0pt, outer sep=0pt, scale=  1.00] at ( 10.19, 50.47) {1};

\path[draw=drawColor,line width= 0.8pt,line join=round,line cap=round] ( 53.68, 15.58) --
	( 53.68, 17.44);

\path[draw=drawColor,line width= 0.8pt,line join=round,line cap=round] ( 15.58, 53.68) --
	( 17.44, 53.68);

\node[text=drawColor,anchor=base,inner sep=0pt, outer sep=0pt, scale=  1.00] at ( 77.85,  6.57) {$\tau_0$};

\node[text=drawColor,anchor=base,inner sep=0pt, outer sep=0pt, scale=  1.00] at (  9.81, 62.34) {$\sigma_0$};

\path[draw=drawColor,line width= 0.8pt,line join=round,line cap=round] ( 77.85, 15.58) --
	( 77.85, 17.44);

\path[draw=drawColor,line width= 0.8pt,line join=round,line cap=round] ( 15.58, 64.83) --
	( 17.44, 64.83);

\path[draw=drawColor,line width= 0.4pt,line join=round,line cap=round] ( 16.51, 64.83) --
	( 77.85, 64.83);

\path[draw=drawColor,line width= 0.4pt,line join=round,line cap=round] ( 77.85, 16.51) --
	( 77.85, 64.83);

\path[draw=drawColor,line width= 0.4pt,line join=round,line cap=round] (128.03,104.38) --
	( 77.85, 64.83);

\node[text=drawColor,anchor=base,inner sep=0pt, outer sep=0pt, scale=  1.00] at ( 47.18, 38.09) {$A_1$};

\node[text=drawColor,anchor=base,inner sep=0pt, outer sep=0pt, scale=  1.00] at (102.01, 38.09) {$A_2$};

\node[text=drawColor,anchor=base,inner sep=0pt, outer sep=0pt, scale=  1.00] at ( 47.18, 88.27) {$A_3$};
\end{scope}
\end{tikzpicture}
 		\caption{}\label{fig1a}
 	\end{subfigure}
 	\begin{subfigure}[b]{0.33\linewidth}
 		\centering
 		% Created by tikzDevice version 0.12 on 2019-03-19 14:44:34
% !TEX encoding = UTF-8 Unicode
\begin{tikzpicture}[x=1pt,y=1pt]
\definecolor{fillColor}{RGB}{255,255,255}
\path[use as bounding box,fill=fillColor,fill opacity=0.00] (0,0) rectangle (144.54,144.54);
\begin{scope}
\path[clip] (  0.00,  0.00) rectangle (144.54,144.54);
\definecolor{drawColor}{RGB}{0,0,0}

\path[draw=drawColor,line width= 0.8pt,line join=round,line cap=round] ( 16.51, 16.51) -- (128.03, 16.51);

\path[draw=drawColor,line width= 0.8pt,line join=round,line cap=round] (124.64, 15.27) --
	(128.03, 16.51) --
	(124.64, 17.74);

\path[draw=drawColor,line width= 0.8pt,line join=round,line cap=round] ( 16.51, 16.51) -- ( 16.51,128.03);

\path[draw=drawColor,line width= 0.8pt,line join=round,line cap=round] ( 17.74,124.64) --
	( 16.51,128.03) --
	( 15.27,124.64);

\node[text=drawColor,anchor=base,inner sep=0pt, outer sep=0pt, scale=  1.00] at (131.38,  9.55) {$\tau$};

\node[text=drawColor,anchor=base,inner sep=0pt, outer sep=0pt, scale=  1.00] at ( 13.16,128.88) {$\sigma$};

\node[text=drawColor,anchor=base,inner sep=0pt, outer sep=0pt, scale=  1.00] at ( 53.68,  5.86) {1};

\node[text=drawColor,anchor=base,inner sep=0pt, outer sep=0pt, scale=  1.00] at ( 10.19, 50.47) {1};

\path[draw=drawColor,line width= 0.8pt,line join=round,line cap=round] ( 53.68, 15.58) --
	( 53.68, 17.44);

\path[draw=drawColor,line width= 0.8pt,line join=round,line cap=round] ( 15.58, 53.68) --
	( 17.44, 53.68);

\node[text=drawColor,anchor=base,inner sep=0pt, outer sep=0pt, scale=  1.00] at ( 44.39,  6.57) {$\tau_1$};

\node[text=drawColor,anchor=base,inner sep=0pt, outer sep=0pt, scale=  1.00] at (  9.81, 34.45) {$\sigma_1$};

\path[draw=drawColor,line width= 0.8pt,line join=round,line cap=round] ( 44.39, 15.58) --
	( 44.39, 17.44);

\path[draw=drawColor,line width= 0.8pt,line join=round,line cap=round] ( 15.58, 36.95) --
	( 17.44, 36.95);

\path[draw=drawColor,line width= 0.4pt,line join=round,line cap=round] ( 16.51, 16.51) --
	( 44.39, 36.95);

\path[draw=drawColor,line width= 0.4pt,line join=round,line cap=round] ( 44.39, 36.95) --
	(128.03, 36.95);

\path[draw=drawColor,line width= 0.4pt,line join=round,line cap=round] ( 44.39, 36.95) --
	( 44.39,128.03);

\path[draw=drawColor,line width= 0.4pt,dash pattern=on 1pt off 3pt ,line join=round,line cap=round] ( 44.39, 36.95) --
	( 16.51, 36.95);

\path[draw=drawColor,line width= 0.4pt,dash pattern=on 1pt off 3pt ,line join=round,line cap=round] ( 44.39, 36.95) --
	( 44.39, 16.51);

\node[text=drawColor,anchor=base,inner sep=0pt, outer sep=0pt, scale=  1.00] at ( 72.27, 69.77) {$\tilde A_1$};

\node[text=drawColor,anchor=base,inner sep=0pt, outer sep=0pt, scale=  1.00] at ( 72.27, 24.23) {$\tilde A_3$};

\node[text=drawColor,anchor=base,inner sep=0pt, outer sep=0pt, scale=  1.00] at ( 30.45, 69.77) {$\tilde A_2$};
\end{scope}
\end{tikzpicture}
 		\caption{}\label{fig1b}
 	\end{subfigure}
 	\begin{subfigure}[b]{0.33\linewidth}
 		\centering
 		% Created by tikzDevice version 0.12 on 2019-03-19 14:48:08
% !TEX encoding = UTF-8 Unicode
\begin{tikzpicture}[x=1pt,y=1pt]
\definecolor{fillColor}{RGB}{255,255,255}
\path[use as bounding box,fill=fillColor,fill opacity=0.00] (0,0) rectangle (144.54,144.54);
\begin{scope}
\path[clip] (  0.00,  0.00) rectangle (144.54,144.54);
\definecolor{drawColor}{RGB}{0,0,0}

\path[draw=drawColor,line width= 0.8pt,line join=round,line cap=round] ( 16.51, 16.51) -- (128.03, 16.51);

\path[draw=drawColor,line width= 0.8pt,line join=round,line cap=round] (124.64, 15.27) --
	(128.03, 16.51) --
	(124.64, 17.74);

\path[draw=drawColor,line width= 0.8pt,line join=round,line cap=round] ( 16.51, 16.51) -- ( 16.51,128.03);

\path[draw=drawColor,line width= 0.8pt,line join=round,line cap=round] ( 17.74,124.64) --
	( 16.51,128.03) --
	( 15.27,124.64);

\node[text=drawColor,anchor=base,inner sep=0pt, outer sep=0pt, scale=  1.00] at (131.38,  9.55) {$\tau$};

\node[text=drawColor,anchor=base,inner sep=0pt, outer sep=0pt, scale=  1.00] at ( 13.16,128.88) {$\sigma$};

\node[text=drawColor,anchor=base,inner sep=0pt, outer sep=0pt, scale=  1.00] at ( 53.68,  5.86) {1};

\node[text=drawColor,anchor=base,inner sep=0pt, outer sep=0pt, scale=  1.00] at ( 10.19, 50.47) {1};

\path[draw=drawColor,line width= 0.8pt,line join=round,line cap=round] ( 53.68, 15.58) --
	( 53.68, 17.44);

\path[draw=drawColor,line width= 0.8pt,line join=round,line cap=round] ( 15.58, 53.68) --
	( 17.44, 53.68);

\path[draw=drawColor,line width= 0.4pt,line join=round,line cap=round] ( 16.51, 16.51) --
	( 53.68, 53.68);

\path[draw=drawColor,line width= 0.4pt,line join=round,line cap=round] ( 53.68, 53.68) --
	(128.03, 53.68);

\path[draw=drawColor,line width= 0.4pt,line join=round,line cap=round] ( 53.68, 53.68) --
	( 53.68,128.03);

\path[draw=drawColor,line width= 0.4pt,dash pattern=on 1pt off 3pt ,line join=round,line cap=round] ( 53.68, 53.68) --
	( 16.51, 53.68);

\path[draw=drawColor,line width= 0.4pt,dash pattern=on 1pt off 3pt ,line join=round,line cap=round] ( 53.68, 53.68) --
	( 53.68, 16.51);

\node[text=drawColor,anchor=base,inner sep=0pt, outer sep=0pt, scale=  1.00] at ( 72.27, 69.77) {$ A_{4,1}$};

\node[text=drawColor,anchor=base,inner sep=0pt, outer sep=0pt, scale=  1.00] at ( 72.27, 32.59) {$ A_{4,3}$};

\node[text=drawColor,anchor=base,inner sep=0pt, outer sep=0pt, scale=  1.00] at ( 35.09, 69.77) {$ A_{4,2}$};
\end{scope}
\end{tikzpicture}
 		\caption{}\label{fig1c}
 	\end{subfigure}
 	\caption{}
 \end{figure}

We have $ \s_0>\t_0>1, \ \  a>1,$ and
 we define the following sets (see Figure \ref{fig1a}) in the first quadrant of the plane $(\t, \s)$:
\begin{eqnarray}\label{sets1}
A_1 &=&\left \{(\t, \s)\in (0, \infty)^2: \t<\t_0, \ \s<\s_0 \right \}, \\  \nonumber
A_2 &=&\left \{(\t, \s)\in (0, \infty)^2: \t>\t_0, \ \s<a\t \right \}, \\  \nonumber
 A_3 &=&\left \{(\t, \s)\in (0, \infty)^2: \s>\s_0, \ \s>a\t  \right \}.
\end{eqnarray}
%This division of the set $(0, \infty)^2$ has very simple structure, it is defined by the point $(\t_0, \s_0)$
%and three lines going through this point - vertical $\{\t=\t_0, \s <\s_0\}$, horizontal $\{\t<\t_0, \s =\s_0\}$, and the line through points $(\t_0, \s_0)$ and $(0, 0)$.
It is easy to see that the function $K_p$ is equal to $1$ on the set $A_p,$ while on the borders between two sets two corresponding functions are equal to one, for example, on the interval $\{ \t=\t_0, \ \s<\s_0\}$ we have  $K_1(\t,\s)=K_2(\t,\s)=1$. At the point $(\t_0, \s_0)$ all three functions $K_p$ are equal to $1$.

Let us denote by ${\tilde S}_{n}(\bt)$ the sum $S_{\bn}(\bt)$ with $\bn=(n, n^\t, n^\s)$ and similarly ${\tilde f}_{n}^{(4)}(\bu)$ and ${\tilde J}_{n}$. We can rewrite (\ref{fbn}) as
\begin{equation}\label{fn1}
\frac{{\tilde f}_{n}^{(4)}(\bu)}{M_{n}^\al(\t, \s)}=\abs{\sum_{l=1}^{d}x_l\left( \sum_{p=1}^3B_p(\bn, \bu, l)\frac{z_{\g_p,n}^{(p)}}{M_{n}(\t, \s)}\frac{ U_{t^{(l)}_p,\g_p}(\sv{n_pu_p},n_p)}{z_{\g_p,n}^{(p)}}  \right) }^\al.
\end{equation}
Having point-wise convergence (\ref{lim5}), in order to apply the Lebesgue dominated convergence theorem we must find majorizing function, but since this is done exactly in the same way as in Example 1, we skip this step. Thus we get
\begin{equation}\label{Jn1}
\frac{{\tilde J}_{n}(\bu)}{n^{1+\t+\s}M_{n}^\al(\t, \s)} \to \int_{-\infty}^{\infty}\int_{-\infty}^{\infty}\int_{-\infty}^{\infty}f(\bu, \t, \s)\dd u_1\dd u_2\dd u_3,
\end{equation}
where
$$
f(\bu, \t, \s)=\abs{\sum_{l=1}^{d}x_l\left( \sum_{p=1}^3 B_p(\bu, l)K_p(\t,\s)H_{\g_p}(u_p,t^{(l)}_p)  \right) }^\al,
$$
$$
B_p(\bu, l)=\ind{\{0\leq u_r \leq t^{(l)}_r, \ 0\leq u_s \leq t^{(l)}_s, r\ne s, r, s \in Q_p\}}.
$$

\begin{prop}\label{prop6} Suppose that we have the sum ${\tilde S}_{n}(\bt)$ of a linear r.f. with the filter (\ref{ci3d}). If $1/\al<\g_1<\g_2<\g_3<1$,  then, as $n \to \infty$,
\begin{equation}\label{lim6}
(n^{1+\t+\s})^{-1/\al}M_{n}^{-1}{\tilde S}_{n}(\bt)\fdd W(\bt, \t, \s)
\end{equation}
where
\begin{equation}\label{lim7}
W(\bt, \t, \s)=\int_{-\infty}^{\infty}\int_{-\infty}^{\infty}\int_{-\infty}^{\infty}\left(\sum_{p=1}^3 B_p(\bu, \bt)K_p(\t,\s)H_{\g_p}(u_p,t_p)   \right)M(\dd u_1,\dd u_2, \dd u_3)
\end{equation}
and
$$
B_p(\bu, \bt)=\ind{\{0\leq u_r \leq t_r, \ 0\leq u_s \leq t_s, r\ne s, r, s \in Q_p\}}.
$$
\end{prop}

Due to the scale transition  we have three types of  limit processes: on each set $A_i$ in the expression (\ref{lim7}) only $K_i(\t, \s)=1$, while two other functions $K_j, j\in Q_i$ are equal to zero;
 on boundaries between any  two of these sets corresponding two functions are equal to $1$; finally, all three functions are equal at the point $(\t_0, \s_0)$: $K_i(\t_0, \s_0)=1, \ i=1, 2, 3$.
We see that in the  case $d=3$ it is difficult to use the terms of well-balanced and unbalanced scaling limits introduced in Definition \ref{sctran}, therefore in \cite{Surg} there was proposed three terms for limit r.f.: {\it well-balanced, partially unbalanced}, and {\it completely unbalanced}. In our context, using this terminology, we have well-balanced limit at the point $(\t_0, \s_0)$, partially unbalanced limits on   boundaries between any  two of  sets $A_i, \ i=1, 2, 3$, and completely unbalanced limits on sets $A_i, \ i=1, 2, 3$.

In this example we took all parameters $\g_i <1, \ i=1, 2, 3,$ which means that we consider long-range dependence along all three axes. We still have the scale transition effect if along one axis we have short-range dependence. Let us consider the case $1/\al<\g_1<\g_2<1<\g_3$, then the normalization quantities $z_{\g_p,n}^{(p)}, \ p=1, 2,$ remain the same, while $z_{\g_3,n}^{(3)}\equiv 1$. Therefore, we get $M_{n}(\t, \s)=\max\{n^{1-\g_1}, n^{\tau (1-\g_2)}\}$, and defining the sets
\begin{eqnarray}\label{set3}
A_4 &=&\left \{(\t, \s)\in (0, \infty)^2: \t<\t_0, \right \}, \\  \nonumber
A_5 &=&\left \{(\t, \s)\in (0, \infty)^2: \t>\t_0,  \right \}, \\  \nonumber
\end{eqnarray}
we get $K_3(\t,\s)\equiv 0, K_1(\t,\s)=1$ on the set $A_4$ and $K_2(\t,\s)=1$ on the set $A_5$. On the half-line (the border between $A_4$ and $A_5$) $\t=\t_0, \s>0$ both functions $K_1$ and $K_2$ equal to $1$.

Similarly, in the case $1/\al<\g_1<\g_3<1<\g_2$ we should get the sets
\begin{eqnarray}\label{set4}
A_6 &=&\left \{(\t, \s)\in (0, \infty)^2: \s<\s_0, \right \}, \\  \nonumber
A_7 &=&\left \{(\t, \s)\in (0, \infty)^2: \s>\s_0,  \right \}, \\  \nonumber
\end{eqnarray}
and $K_2(\t,\s)\equiv 0, K_1(\t,\s)=1$ on the set $A_6$ and $K_3(\t,\s)=1$ on the set $A_7$.

\bigskip

\subsection{Comparison of Example 4 with results from \cite{Surg}}

We can compare our Example 4 with the results from \cite{Surg}, and to this aim we formulate main results from \cite{Surg}, using the notation, most close to our notation. In \cite{Surg} a linear r.f.  $$
X_\bk=\sum_{\bi\in\ZZ^3}c_{\bk-\bi}\e_{\bi}, \ \bk \in \ZZ^3,
$$
with a filter
\begin{equation}\label{DSfilter}
c_{\bi}=\frac{g(\bi)}{(\sum_{j=1}^3a_j|i_j|_+^{\g_j/\nu})^\nu}, \ \bi \in \ZZ^3,
\end{equation}
is considered. Here $\e_{\bi}, \ \bi \in \ZZ^3$ are i.i.d. random variables with mean zero and unit variance (this would correspond to our case $\al=2$),  $|a|_+=\max (|a|, 1)$, $g(\bi),  \ \bi \in \ZZ^3,$ is bounded function and $\lim_{|\bi|\to \infty} g(\bi):=g_\infty \in (0. \infty)$, $a_j>0, j=1, 2, 3, \ \nu>0,$ and parameters $\g_j$ satisfy the following condition
$$
1<Q:=\sum_{j=1}^3\frac{1}{\g_j}<2.
$$
This condition guarantees that
$$
\sum_{\bi\in\ZZ^3}|c_{\bi}|^2<\infty \quad  {\rm and} \quad  \sum_{\bi\in\ZZ^3}|c_{\bi}|=\infty,
$$
i.e., r.f. is a stationary with finite variance and with long-range dependence. Let us note that coefficients of the filter along axes decay at the same rates as in our Example 4, namely, if we denote $\bi^{(j)}$ vector $\bi$ with $i_k=0$ for $k\ne j$, then it is easy to see that $c_{\bi^{(j)}}=O(|i_j|_+^{-\g_j})$. If it would be possible to take the function $g$ equal to zero at all points of $\ZZ^3$ except axes, then  our Example 4 would follow from the results in \cite{Surg}. But this is not possible due to the requirement that the limits of  $g(\bt)$, as $|\bt|\to \infty$ must be the same and positive. As a matter of fact, in the proofs in \cite{Surg} it is assumed without loss of generality that $g(\bt)\equiv 1$. The values of the r.f. $X_\bk$ are summed over rectangles as in (\ref{field}) with $n_i=n^{q_i}, \ i=1, 2, 3,$ but essential parameters are the ratios between $q_i$, therefore, in order to get the same notation as our one, we can set $q_1=1, q_2=\t, q_3=\s$. With this notation the balance conditions in \cite{Surg} are expressed by means of $\t_0=\g_1/\g_2$ and $\s_0=\g_1/\g_3$, and in Fig 1 in \cite{Surg} there is given the partition of the quadrant $\{\t>0, \s>0 \}$ by means of balance conditions. This picture has the same structure as our Fig 1, if we take all sets $A_1 - A_7$. But this is the only similarity between results in \cite{Surg} and our Example 4. The regions of parameters $\g_i$, the values of $\t_0$ and $\s_0$, and limit  r.f.'s in \cite{Surg} and Example 4 are different. This can be explained by the fact that in \cite{Surg} the filter of a r.f. under consideration is "three-dimensional" (in the sense that coefficients  are non-zero over all $\ZZ^3$) while in our examples filters are "one-dimensional" (non-zero only on axes or some lines). Such simple filters allow to understand better the mechanism of scaling transition and motivated more general definitions of scale transition comparing with original definition given in \cite{Puplinskaite2}.

\subsection{Example 5}

 In Example 2 (the case $d=2$) we had shown, that long-range dependence connected with the requirement that both exponents $\g_i<1,  i=1, 2,$ is not necessary condition for the scaling transition, and the scale transition can be observed in the case of  negative dependence. The same can be shown in the case $d=3$ and to do this
  one needs to make small change in the filter (\ref{ci3d}). Namely,
 we take the same three sequences of positive numbers $a_i^{(j)},  \ i\ge 1, \ j=1, 2,3,$ but now we assume $\max (1, 1/\al)<\g_1<\g_2<\g_3<1+1/\al$ and we define $c_{(0,0,0)}=a_1^{(0)}+a_2^{(0)}+a_3^{(0)}$ with some numbers $a_j^{(0)}$ in such a way that  the following conditions
\begin{equation}\label{sum3dZero}
\sum_{\bi\ge \B0}c_{\bi}=0
\end{equation}
and $\sum_{i=0}^\infty a_j(i)=0, \ j=1, 2, 3,$ are satisfied.  It is possible to show that in this case we get very similar picture as in Example 4, but since the proofs are very similar to those used in Examples 2 and 4, we shall give only the final result.

 Now the point $(\t_1, \s_1),$ which determines the scaling transition, will be
$$
 \t_1=\frac{ \g_1-1}{ \g_2-1}, \ \s_1=\frac{\g_1-1}{\g_3-1}, \ a_1=\frac{ \g_2-1}{\g_3-1} %\ b=\frac{ 1-\g_1}{ 1-\g_3},
 $$
and these quantities satisfy $0< \s_1<\t_1<1, \ \  a_1<1.$ As in Example 4 we get three sets (see Figure \ref{fig1b}) in the first quadrant of the plane $(\t, \s)$, in which we have different limit fields:
\begin{eqnarray}\label{sets2}
{\tilde A}_1 &=&\left \{(\t, \s)\in (0, \infty)^2: \t>\t_0, \ \s>\s_0 \right \}, \\  \nonumber
{\tilde A}_2 &=&\left \{(\t, \s)\in (0, \infty)^2: \t<\t_0, \ \s>a\t \right \}, \\  \nonumber
{\tilde A}_3 &=&\left \{(\t, \s)\in (0, \infty)^2: \s<\s_0, \ \s<a\t  \right \}.
\end{eqnarray}

Thus, we get the picture of the scale transition which is in a sense inverse to the picture of Example 4, and  the structure of limit r.f. in this case is very similar to (\ref{lim7}).

 \subsection{Example 6}

 Considering the case $d=2$ (Examples 1-3)  we saw that filters with coefficients or on axes either on diagonal give us one point of the scaling transition, while combining diagonal with axis  we got two points of the scaling transition. Therefore, in the case $d=3$, in order to get more complicated picture of the scaling transition, it is natural  to consider the linear r.f. (\ref{field3d}) with the following filter
\begin{equation}\label{ci3ddiag}
c_{\bi}=
\begin{cases}
a_{1}(i_1), \text{ if } i_1\geq 1, i_2=i_3=0,\\
a_{2}(i_2), \text{ if } i_2\geq 1, i_1=i_3=0,\\
a_{3}(i_3), \text{ if } i_3\geq 1, i_2=i_1=0,\\
a_{4}(i), \text{ if } ,  i_3=i_2=i_1=i\ge 1, \\
0, \text{ elsewhere}.
\end{cases}
\end{equation}
where sequences $a_{p}(i), p=1, 2, 3,$  are as in Example 4 and the fourth sequence is $a_{4}(0)=0, \ a_{4}(i)=(1+i)^{-\g_4}, \ 1/\al<\g_4 <1$. If it was possible to compare Example 4 with results for linear r.f. with filter (\ref{DSfilter}), since coefficients of this filter have different decay rates along axes, it is easy to see that coefficients of this filter on diagonal decay with the same rate as coefficients on the axis with the slowest rate of decay, therefore Example 6 cannot be compared with results from \cite{Surg}.
For a moment we do not relate the parameter $\g_4$ with parameters $\g_i, \ i=1, 2, 3.$ Since the filter of this example is obtained by combining filter from Examples 4 and adding coefficients on the diagonal, as in Example 3 ( in the case $d=2$) it is easy to write the ch.f. for $A_{\bn}^{-1}\sum_{l=1}^{d}x_lS_{\bn}(\bt^{(l)})$ and to get formulae analogous to (\ref{Jbn}), (\ref{fbn}). Thus we need to investigate the quantity $I_{\bn}:=A_{\bn}^{-\al}J_{\bn},$ where
\begin{equation}\label{Jbn2}
J_{\bn}=n_1n_2n_3
\int_{-\infty}^{\infty}\int_{-\infty}^{\infty}\int_{-\infty}^{\infty}f_{\bn}^{(4)}(\bu)\dd u_1\dd u_2\dd u_3,
\end{equation}
and
\begin{equation}\label{fbn2}
f_{\bn}^{(5)}(\bu)=\abs{\sum_{l=1}^{d}x_l\left( \sum_{p=1}^3 B_p(\bn, \bu, l) U_{t^{(l)}_p,\g_p}(\sv{n_pu_p},n_p)+\prod_{j=1}^3\ind{\{0\leq \sv{n_ju_j}\leq \sv{n_j t^{(l)}_j}\}}Z_{\bt^{(l)}, \g_4}(\sv{\bn\bu}, \bn)  \right) }^\al,
\end{equation}
\begin{equation}\label{zdiag3}
Z_{\bt^{(l)}, \g_4}(\bi, \bn):=\sum_{k=\max\left( 0,-i_1,-i_2, -i_3 \right)}^{\min\left( \sv{n_jt^{(l)}_j}-i_j, j=1, 2, 3, \right)}a_{4}(k).
\end{equation}
The expression of $f_{\bn}^{(4)}(\bu)$ is similar to that given in (\ref{fbn}) (differs by one additional term), quantities, present in (\ref{fbn2}) are defined in (\ref{Bpbn}), (\ref{Sngamma}).
As in examples above, we set $n_1=n, \ n_2=n^\t,\ n_3=n^\s$ and $1/\al<\g_1<\g_2<\g_3<1$. Now we must find the right normalization for four terms present in (\ref{fbn2}), but for three terms with $U_{t^{(l)}_p,\g_p}(\sv{n_pu_p},n_p)$ normalization is obtained in Example 4 and is given by quantities $z_{\g_p,n}^{(p)}$ (see (\ref{lim5})) and by sets $A_i,  i=1, 2, 3,$ (see (\ref{sets1})). Namely, from Example 4 we have that for $(\t, \s)\in A_p$ the normalization for the first three terms in  (\ref{fbn2}) is $z_{\g_p,n}^{(p)}=n^{s_p}, \ p=1, 2, 3,$ where
\begin{equation}\label{normalization1}
s_1=1-\g_1, \ s_2=\t(1-\g_2), \    s_3=\s(1-\g_3).
\end{equation}
Now let us consider normalization for the fourth term in (\ref{fbn2}).
It is easy to see that the normalization for $Z_{\bt^{(l)}, \g_4}(\bi, \bn)$ depends on the relations between parameters $\t$ and $\s$, since for the investigation of the growth of this sum we must compare quantities $n_i, \ i=1, 2, 3.$ We shall skip the procedure of this comparison and provide the final result.
 As in Example 4, the growth of the sum in (\ref{zdiag3}) is different on three sets of the possible values of $(\t, \s)$, and    the  division of the area of possible values of parameters $(\t, \s) \in (0, \infty)^2$  into three sets (see Figure \ref{fig1c}), is by  the point  $(\t_1, \s_1)$, with   $\t_1=\s_1=1$:
 \begin{eqnarray}\label{setslyg}
A_{4, 1} &=& \{(\t, \s)\in (0, \infty)^2: \t>1, \ \s>1 \}, \\ \nonumber
A_{4, 2} &=& \{(\t, \s)\in (0, \infty)^2: \s>\t, \ \t<1  \}, \\  \nonumber
A_{4, 3} &=& \{(\t, \s)\in (0, \infty)^2: \s<\t, \ \s<1 \}.
\end{eqnarray}
Normalization for (the growth of)  $Z_{\bt^{(l)}, \g_4}(\bi, \bn)$,  according to our assumption, is a function of $n$, and on the set  $A_{4, i}$,  is $z_{\g_4,n}^{(4, i)}=n^{s_{4, i}}, i=1, 2, 3,$  where
\begin{equation}\label{normalization2}
s_{4, 1}=1-\g_4,
s_{4, 2}=\t(1-\g_4),
s_{4, 3}=\s(1-\g_4).
\end{equation}

 It is clear that, in order to find normalization for the function (\ref{fbn2}), we must find $M(n; \t, \s)=\max \{z_{\g_p,n}^{(p)}, \ p=1, 2, 3, z_{\g_4,n}^{(4, i)}, \ i=1, 2, 3,\}$. For this aim we consider the intersections of sets $A_p$ and $A_{4, j}, $ with all possible combinations of $1\le p, j\le 3$. Since the point $(\t_0, \s_0)$ is defined by parameters $\g_i$ (and we  fixed the relation $1/\al<\g_1<\g_2<\g_3<1$) and  $(\t_1, \s_1)$ is independent of all $\g_i$, we have  seven sets (two intersections $A_2\cap A_{4, 2}, \ A_3\cap A_{4, 3}$ are empty, since $\t_0>\t_1=1, \ \s_0>\s_1=1$):
$$
B_1=A_2\cap A_{4, 1}, \ B_2=A_3\cap A_{4, 1}, \ B_3=A_3\cap A_{4, 2}, \ B_4=A_1\cap A_{4, 2}
$$
$$
 B_5=A_1\cap A_{4, 1}, \ B_6=A_1\cap A_{4, 3}, \ B_7=A_2\cap A_{4, 3}
$$
Now in each set $B_i$ we must compare only two exponents - if $ B_i=A_p\cap A_{4, j}$, we must compare $s_p$ and $ s_{4, j}$.  Till now the choice of the parameter $\g_4$ was arbitrary in the interval $(1/\al, 1)$, in which all three $\g_i,  i=1, 2, 3,$ are located. It turns out that the position, where $\g_4$ is located, is important. Let us consider the case $1/\al<\g_1<\g_2<\g_3<\g_4<1$, i.e., the sequence on the diagonal is decreasing most rapidly. Then easy calculations show that in this case in all sets $B_i$ the prevailing (i.e., bigger) exponents are $s_j$: in the sets $B_4, B_5, B_6$ the prevailing is $s_1$,
 in $B_1, B_7$  - $s_2$ and in $B_2$ and $B_3$ the  exponent $s_3$ is bigger than $s_{4, 1}$ and $s_{4, 2}$, respectively. Since $A_1=B_4\cup B_5\cup B_6, \ A_2=B_1 \cup B_7$ and $A_3=B_2\cup B_3$, we see that the scaling transition in this case, where the sequence on the diagonal, comparing with sequences on axes, is decreasing most rapidly, is completely determined by the sets $A_i, \ i=1, 2, 3,$ and the filter coefficients on axes.

 As some surprise for us it was the fact that the same picture of scaling transition, i.e., defined by sets  $A_i$ and exponents $s_i$ \ $i=1, 2, 3,$  we get also in two cases $1/\al<\g_1<\g_2<\g_4<\g_3<1$ and $1/\al<\g_1<\g_4<\g_2<\g_3<1$. This is obtained in the same way - comparing exponents of normalizing constants on each of the sets $B_i, \ 1\le i \le 7$.

 Let us consider  the last  case, where the sequence on the diagonal is decreasing most slowly, i.e., $1/\al <\g_4<\g_1<\g_2<\g_3<1$. This case gives us the  picture, which we are looking for:  the scaling transition is defined by all sets $A_i$ and $A_{4, i}, \ i=1, 2, 3$ and even some sets $B_j$ are divided by new lines. We provide the complete analysis of this case. We recall that we had introduced notations $\t_0=(1-\g_1)/(1-\g_2), \s_0=(1-\g_1)/(1-\g_3), \ \t_1=\s_1=1,$ additionally we  denote $\t_2=(1-\g_4)/(1-\g_2), \ \t_3=(1-\g_1)/(1-\g_4), \ \s_2=(1-\g_4)/(1-\g_3). $ We shall use without special mentioning  the inequalities
 $$
 1-\frac{1}{\al}>1-\g_4>1-\g_1>1-\g_2>1-\g_3>0.
 $$
Let us consider the set $B_1.$ We must compare $\t(1-\g_2)$   and $1-\g_4$. In the set $B_1$ we have $\t>\t_0>1$ and since $1-\g_4>1-\g_1$, therefore, for $\t_0<\t<\t_2$, we get that $1-\g_4>\t(1-\g_2)$,  while if  $\t>\t_2$, then $1-\g_4<\t(1-\g_2)$. Thus, we get that the set $B_1$ is divided into two parts by the vertical line going through the point $(\t_2, 0)$.

 In the set $B_2$  we must compare $\s(1-\g_3)$   and $1-\g_4$. In this set  we have $\s>\s_0>1$ and, again using  $1-\g_4>1-\g_1$, we get that $1-\g_4>\s(1-\g_3)$, for $\s_0<\t<\s_2$, while if  $\s>\s_2$, then $1-\g_4<\s(1-\g_3)$. We get that the set $B_2$ is divided into two parts by the horizontal line going through the point $(0, \s_2).$

In the set $B_3$  we must compare $\s(1-\g_3)$   and $\t(1-\g_4)$. In this set  we have $\s>\t(1-\g_2)/(1-\g_3)$ and since   $1-\g_4>1-\g_2$ we get that $\t(1-\g_4)>\s(1-\g_3)$, for $\t(1-\g_2)/(1-\g_3)<\s<\s_2\t$, while if  $\s>\s_2\t$, then $\t(1-\g_4)<\s(1-\g_3)$. We get that the set $B_3$ is divided into two parts by the  line  $\s=\s_2 \t).$

In the set $B_4$  we must compare $1-\g_1$   and $\t(1-\g_4)$. It is easy to see that if $\t_3<\t<1,$ then $\t(1-\g_4)>1-\g_1$, while for $\t_3>\t$ the opposite inequality holds, this means that the set $B_4$ is divided into two parts by the vertical line going through the point $(\t_3, 0)$. Also from our assumption easily follows that in the set $B_5$ we have $1-\g_4>1-\g_1$

In the set $B_6$  we must compare $\s(1-\g_4)$   and $1-\g_1$. It is not difficult to verify that this set is divided into two parts by the horizontal line going through the point $(0, \s_3),$ namely, if $\s_3<\s<1$ then $\s(1-\g_4)>1-\g_1$, while if $\s_3>\s$ then $\s(1-\g_4)<1-\g_1$.

Finally, in the set $B_7$,   comparing $\s(1-\g_4)$   and $\t(1-\g_2)$, we get that this set is divided into two parts by the  line $\s=\t\t_2^{-1}$, and if $\t\t_2^{-1}<\s<1$ then $\s(1-\g_4)>\t(1-\g_2)$, while if $\t\t_2^{-1}>\s$ then $\s(1-\g_4)<\t(1-\g_2)$.

\begin{figure}[!ht]
	\centering
	% Created by tikzDevice version 0.12 on 2019-03-21 07:30:02
% !TEX encoding = UTF-8 Unicode
\begin{tikzpicture}[x=1pt,y=1pt]
\definecolor{fillColor}{RGB}{255,255,255}
\path[use as bounding box,fill=fillColor,fill opacity=0.00] (0,0) rectangle (289.08,289.08);
\begin{scope}
\path[clip] (  0.00,  0.00) rectangle (289.08,289.08);
\definecolor{drawColor}{RGB}{0,0,0}

\path[draw=drawColor,line width= 0.8pt,line join=round,line cap=round] ( 22.87, 22.87) -- (266.21, 22.87);

\path[draw=drawColor,line width= 0.8pt,line join=round,line cap=round] (262.81, 21.64) --
	(266.21, 22.87) --
	(262.81, 24.11);

\path[draw=drawColor,line width= 0.8pt,line join=round,line cap=round] ( 22.87, 22.87) -- ( 22.87,266.21);

\path[draw=drawColor,line width= 0.8pt,line join=round,line cap=round] ( 24.11,262.81) --
	( 22.87,266.21) --
	( 21.64,262.81);

\node[text=drawColor,anchor=base,inner sep=0pt, outer sep=0pt, scale=  1.00] at (269.86, 15.51) {$\tau$};

\node[text=drawColor,anchor=base,inner sep=0pt, outer sep=0pt, scale=  1.00] at ( 19.22,267.36) {$\sigma$};

\node[text=drawColor,anchor=base,inner sep=0pt, outer sep=0pt, scale=  1.00] at (103.98, 11.55) {1};

\node[text=drawColor,anchor=base,inner sep=0pt, outer sep=0pt, scale=  1.00] at ( 15.98,100.78) {1};

\path[draw=drawColor,line width= 0.8pt,line join=round,line cap=round] (103.98, 21.86) --
	(103.98, 23.89);

\path[draw=drawColor,line width= 0.8pt,line join=round,line cap=round] ( 21.86,103.98) --
	( 23.89,103.98);

\path[draw=drawColor,line width= 0.4pt,dash pattern=on 1pt off 3pt ,line join=round,line cap=round] ( 83.71, 22.87) --
	( 83.71, 83.71);

\path[draw=drawColor,line width= 0.4pt,dash pattern=on 1pt off 3pt ,line join=round,line cap=round] ( 22.87, 83.71) --
	( 83.71, 83.71);

\path[draw=drawColor,line width= 0.4pt,dash pattern=on 1pt off 3pt ,line join=round,line cap=round] (103.98,103.98) --
	( 22.87,103.98);

\path[draw=drawColor,line width= 0.4pt,dash pattern=on 1pt off 3pt ,line join=round,line cap=round] (103.98,103.98) --
	(103.98, 22.87);

\path[draw=drawColor,line width= 0.4pt,dash pattern=on 1pt off 3pt ,line join=round,line cap=round] (103.98,163.47) --
	( 22.87,163.47);

\path[draw=drawColor,line width= 0.4pt,dash pattern=on 1pt off 3pt ,line join=round,line cap=round] (201.32,163.47) --
	(201.32, 22.87);
\definecolor{drawColor}{RGB}{190,190,190}

\path[draw=drawColor,line width= 0.4pt,dash pattern=on 4pt off 4pt ,line join=round,line cap=round] ( 22.87, 22.87) --
	( 83.71, 83.71);

\path[draw=drawColor,line width= 0.4pt,dash pattern=on 4pt off 4pt ,line join=round,line cap=round] (103.98,163.47) --
	(103.98,266.21);

\path[draw=drawColor,line width= 0.4pt,dash pattern=on 4pt off 4pt ,line join=round,line cap=round] (201.32,103.98) --
	(266.21,103.98);

\path[draw=drawColor,line width= 0.4pt,dash pattern=on 1pt off 3pt on 4pt off 3pt ,line join=round,line cap=round] ( 83.71,128.32) --
	(156.71,128.32);

\path[draw=drawColor,line width= 0.4pt,dash pattern=on 1pt off 3pt on 4pt off 3pt ,line join=round,line cap=round] (156.71, 83.71) --
	(156.71,128.32);

\path[draw=drawColor,line width= 0.4pt,dash pattern=on 1pt off 3pt on 4pt off 3pt ,line join=round,line cap=round] (201.32,163.47) --
	(156.71,128.32);
\definecolor{drawColor}{RGB}{255,0,0}

\path[draw=drawColor,line width= 0.4pt,line join=round,line cap=round] (103.98,103.98) --
	( 83.71, 83.71);

\path[draw=drawColor,line width= 0.4pt,line join=round,line cap=round] (156.71, 83.71) --
	( 83.71, 83.71);

\path[draw=drawColor,line width= 0.4pt,line join=round,line cap=round] ( 83.71,128.32) --
	( 83.71, 83.71);

\path[draw=drawColor,line width= 0.4pt,line join=round,line cap=round] (103.98,163.47) --
	(103.98,103.98);

\path[draw=drawColor,line width= 0.4pt,line join=round,line cap=round] (201.32,103.98) --
	(103.98,103.98);

\path[draw=drawColor,line width= 0.4pt,line join=round,line cap=round] (201.32,163.47) --
	(201.32,103.98);

\path[draw=drawColor,line width= 0.4pt,line join=round,line cap=round] (201.32,163.47) --
	(103.98,163.47);

\path[draw=drawColor,line width= 0.4pt,line join=round,line cap=round] (156.71, 83.71) --
	(201.32,103.98);

\path[draw=drawColor,line width= 0.4pt,line join=round,line cap=round] ( 83.71,128.32) --
	(103.98,163.47);

\path[draw=drawColor,line width= 0.4pt,line join=round,line cap=round] (156.71, 83.71) --
	(201.32,103.98);

\path[draw=drawColor,line width= 0.4pt,line join=round,line cap=round] ( 83.71,128.32) --
	(103.98,163.47);

\path[draw=drawColor,line width= 0.4pt,line join=round,line cap=round] (201.32,163.47) --
	(266.21,214.59);

\path[draw=drawColor,line width= 0.4pt,line join=round,line cap=round] (156.71, 22.87) --
	(156.71, 83.71);

\path[draw=drawColor,line width= 0.4pt,line join=round,line cap=round] ( 22.87,128.32) --
	( 83.71,128.32);
\definecolor{drawColor}{RGB}{0,0,0}

\path[draw=drawColor,line width= 0.8pt,line join=round,line cap=round] (156.71, 21.86) --
	(156.71, 23.89);

\path[draw=drawColor,line width= 0.8pt,line join=round,line cap=round] ( 21.86,128.32) --
	( 23.89,128.32);

\node[text=drawColor,anchor=base,inner sep=0pt, outer sep=0pt, scale=  1.00] at (156.71, 12.26) {$\tau_0$};

\node[text=drawColor,anchor=base,inner sep=0pt, outer sep=0pt, scale=  1.00] at ( 15.98,125.82) {$\sigma_0$};

\path[draw=drawColor,line width= 0.8pt,line join=round,line cap=round] ( 83.71, 21.86) --
	( 83.71, 23.89);

\path[draw=drawColor,line width= 0.8pt,line join=round,line cap=round] ( 21.86, 83.71) --
	( 23.89, 83.71);

\node[text=drawColor,anchor=base,inner sep=0pt, outer sep=0pt, scale=  1.00] at ( 83.71, 12.26) {$\tau_3$};

\node[text=drawColor,anchor=base,inner sep=0pt, outer sep=0pt, scale=  1.00] at ( 15.98, 81.21) {$\sigma_3$};

\path[draw=drawColor,line width= 0.8pt,line join=round,line cap=round] (201.32, 21.86) --
	(201.32, 23.89);

\path[draw=drawColor,line width= 0.8pt,line join=round,line cap=round] ( 21.86,163.47) --
	( 23.89,163.47);

\node[text=drawColor,anchor=base,inner sep=0pt, outer sep=0pt, scale=  1.00] at (201.32, 12.26) {$\tau_2$};

\node[text=drawColor,anchor=base,inner sep=0pt, outer sep=0pt, scale=  1.00] at ( 15.98,160.97) {$\sigma_2$};

\path[draw=drawColor,line width= 1.6pt,line join=round,line cap=round] ( 83.71, 83.71) --
	( 83.71, 83.71);

\path[draw=drawColor,line width= 1.6pt,line join=round,line cap=round] ( 83.71,128.32) --
	( 83.71,128.32);

\path[draw=drawColor,line width= 1.6pt,line join=round,line cap=round] (103.98,163.47) --
	(103.98,163.47);

\path[draw=drawColor,line width= 1.6pt,line join=round,line cap=round] (103.98,103.98) --
	(103.98,103.98);

\path[draw=drawColor,line width= 1.6pt,line join=round,line cap=round] (156.71, 83.71) --
	(156.71, 83.71);

\path[draw=drawColor,line width= 1.6pt,line join=round,line cap=round] (201.32,103.98) --
	(201.32,103.98);

\path[draw=drawColor,line width= 1.6pt,line join=round,line cap=round] (201.32,163.47) --
	(201.32,163.47);

\node[text=drawColor,anchor=base,inner sep=0pt, outer sep=0pt, scale=  1.00] at ( 77.62, 86.35) {A};

\node[text=drawColor,anchor=base,inner sep=0pt, outer sep=0pt, scale=  1.00] at ( 77.62,130.96) {B};

\node[text=drawColor,anchor=base,inner sep=0pt, outer sep=0pt, scale=  1.00] at ( 97.90,166.11) {C};

\node[text=drawColor,anchor=base,inner sep=0pt, outer sep=0pt, scale=  1.00] at (110.07,106.62) {G};

\node[text=drawColor,anchor=base,inner sep=0pt, outer sep=0pt, scale=  1.00] at (162.79, 74.18) {F};

\node[text=drawColor,anchor=base,inner sep=0pt, outer sep=0pt, scale=  1.00] at (207.40, 94.46) {E};

\node[text=drawColor,anchor=base,inner sep=0pt, outer sep=0pt, scale=  1.00] at (195.23,166.11) {D};

\node[text=drawColor,anchor=base,inner sep=0pt, outer sep=0pt, scale=  1.00] at (150.62, 25.51) {K};

\node[text=drawColor,anchor=base,inner sep=0pt, outer sep=0pt, scale=  1.00] at ( 28.96,130.96) {H};

\node[text=drawColor,anchor=base,inner sep=0pt, outer sep=0pt, scale=  1.00] at ( 63.43, 60.93) {$ s_{1}$};

\node[text=drawColor,anchor=base,inner sep=0pt, outer sep=0pt, scale=  1.00] at (231.73, 90.67) {$ s_{2}$};

\node[text=drawColor,anchor=base,inner sep=0pt, outer sep=0pt, scale=  1.00] at (112.10,215.04) {$ s_{3}$};

\node[text=drawColor,anchor=base,inner sep=0pt, outer sep=0pt, scale=  1.00] at (142.51,113.65) {$ s_{4,1}$};

\node[text=drawColor,anchor=base,inner sep=0pt, outer sep=0pt, scale=  1.00] at (142.51, 91.35) {$ s_{4,3}$};

\node[text=drawColor,anchor=base,inner sep=0pt, outer sep=0pt, scale=  1.00] at ( 93.85,113.65) {$ s_{4,2}$};
\end{scope}
\end{tikzpicture}
	\caption{}\label{fig4}
\end{figure}

Collecting all these facts we have the following picture (see Figure \ref{fig4}). Let us denote the points $A=(\t_3, \s_3),\  B=(\t_3, \s_0),\ C=(\t_1, \s_2),\ D=(\t_2, \s_2),\ E=(\t_2, \s_0),\ F=(\t_0, \s_3),\ G=(\t_0, \s_0),\ H=(0, \s_0),\ K=(\t_0, 0), \ O=(0, 0),$ then inside the set $ABCDEF$ sets $A_{4, j}$ and exponents $ s_{4, j}$ are dominating: $ s_{4, 1}$ is dominating in $GCDE$, $ s_{4, 2}$   is dominating in $ABCG$ and $ s_{4, 3}$ is dominating in $AGEF$. Outside of the set $ABCDEF$ sets $A_{ j}$ and exponents $ s_{ j}$ are dominating: to the top of the broken line $HBCD$ the exponent $s_3$ is the biggest, to the right from
the broken line $DEFK$ the exponent $s_2$ is dominating and in the set $HBAFKO$ $s_1$ is the biggest.

After this analysis we have the normalizing sequence $A_n=n^{(1+\t+\s)/\al}M(n; \t, \s)$.   Then the final step - finding the limit  distributions -  is carried in the same way as in Example 4, introducing normalization for each term in (\ref{fbn2}) and functions of type (\ref{Kp}). We shall get that inside each set in Figure \ref{fig4} only one function $K_p$ will be equal to $1$, while on lines which serve as border lines we shall have two functions equal to one, while at points $A, B, C, D, E, F, G$ we shall have three functions equal to $1$, since at each of these points three sets with different prevailing exponents meet.

  \subsection{Example 7}

 In the case $d=3$ we consider sums (\ref{sum13d}) indexed with three-dimensional indices $\bn =(n_1, n_2, n_3)$. Clearly, we have the same first possibility as in the case $d=2$, namely to consider limit theorems assuming only $\bn \to \infty$ (we recall that this means that $\min (n_i, \ i=1, 2, 3,) \to \infty$) , but for the second possibility we have more freedom. In Examples 4-6 we assumed that $n_2$ and $n_3$ are functions of $n_1$, more precisely, $n_2=n_1^\t, n_3=n_1^\s$. But in the relation $\min (n_i, \ i=1, 2, 3,) \to \infty$ we may assume that $(n_1, n_2) \to \infty$ and $n_3 = f(n_1, n_2)$ with some function $f: \ZZ_+^2 \to \ZZ_+$ such that $f(n_1, n_2) \to \infty)$, as $(n_1, n_2) \to \infty$. One can hope that for some linear fields and some functions $f$ we can have the scaling transition. The most natural functions to begin with are $[(n_1n_2)^\t], \ [n_1^\t n_2^\s], \  [n_1^\t]+[n_2^\s]$, where $\t$ and $\s$ are positive numbers. We shall take the simple function $f(n_1, n_2)=[n_1^\t n_2^\s]$ and we shall look for a filter of a linear field to get the scaling transition. Since we know  that in the case where filter coefficients are expressed as factors of two sequences ($c_{i, j}=a_ib_j$ in the case $d=2$) there is no scale transition, we can try the following filter. We choose three sequences $a_i(j)=(1+j)^{\g_i},  \ j\ge 1, a_i(0)=0, \ i=1, 2,3,$
 and define
\begin{equation}\label{ci3dfakt}
c_{\bi}=
\begin{cases}
a_1(i_1)a_2(i_2), \text{ if } i_3=0,\\
a_3(i_3), \text{ if } i_1=i_2= 0, i_3\ge 1,\\
0, \text{ elsewhere}.
\end{cases}
\end{equation}
This can be written as
\begin{equation*}
c_{\bi}=a_1(i_1)a_2(i_2)\ind{i_3=0}+a_3(i_3)\ind{i_1=i_2=0},
\end{equation*}
Again we can start from formulae (\ref{generalJn})-(\ref{generalfbn}). Taking into account the expression of the filter (\ref{ci3dfakt}) after some transformations we can obtain
\begin{align}\label{fbn3d fakt}
f_{\bn}^{(6)}\left(\bu,\bt^{(j)}\right) &=\sum_{\bk\in\ZZ^d}
	c_{\bk}\ind{[(-\sv{\bn\bu})\vee \B0\leq\bk\leq\bn\bt^{(j)}-\sv{\bn\bu}]}\\ \nonumber
  &= \sum_{i_1, i_2=0}^\infty a_1(i_1)a_2(i_2)\ind{[(-\sv{n_q u_q})\vee 0\leq k_q\leq n_q t_q^{(j)}-\sv{n_q u_q}, q=1, 2]}\ind{[-\sv{n_3 u_3} \leq 0 \leq n_3 t_3^{(j)}-\sv{n_3 u_3}]}\\ \nonumber
  &+\sum_{i_3=1}^\infty a_3(i_3)\ind{[(-\sv{n_q u_q})\leq 0\leq n_q t_q^{(j)}-\sv{n_q u_q}, q=1, 2]}\ind{[-n_3 u_3\vee 0 \leq 0 \leq n_3 t_3^{(j)}-{n_3 u_3}+n_3]}
\end{align}
Since the first double sum can be written as product of two sums and the second sum is the same as considered in previous examples, using the notation (\ref{Sngamma}), we can write
\begin{align}\label{fbnex8}
f_{\bn}^{(6)}\left(\bu,\bt^{(j)}\right) &=U_{t_1^{(j)},\g_1}(\sv{n_1u_1},n_1)U_{t_2^{(j)},\g_2}(\sv{n_2u_2},n_2)\ind{[0\leq\sv{n_3 u_3}  \leq n_3 t_3^{(j)}]}\\ \nonumber
                                  &+ U_{t_3^{(j)},\g_3}(\sv{n_3u_3},n_3)\ind{[0\leq\sv{n_i u_i}  \leq n_i t_i^{(j)}, \ i=1, 2]}
\end{align}

%\begin{align*}\label{}
%f_{\bn}\left(\bu,\bt^{(j)}\right)& =\sum_{(-\sv{n_1 u_1})\vee 0\leq k_1\leq n_1 t_1^{(j)}-\sv{n_1 u_1}}^\infty a_1(k_1) \sum_{(-\sv{n_2 u_2})\vee 0\leq k_2\leq n_2 t_2^{(j)}-\sv{n_2 u_2}}^\infty a_2(k_2)\ind{[0\leq\sv{n_3 u_3}  \leq n_3 t_3^{(j)}]}\\
% & +\sum_{(-\sv{n_3 u_3})\vee 0\leq k_3\leq n_3 t_3^{(j)}-\sv{n_3 u_3}}^\infty a_3(k_3).
%\end{align*}
From Proposition \ref{prop1} we know that the right normalization for the first term is $z_{\g_1,n_1}z_{\g_2,n_2}$, while $z_{\g_3,n_3}$  is normalization sequence for the second term. Therefore, we must investigate the quantity
\begin{equation}\label{Mbn}
M_{\bn}=M_{\bn}(\g_i, \g_2, \g_3)= \max (z_{\g_1,n_1}z_{\g_2,n_2}, \ z_{\g_3,n_3}),
\end{equation}
and this quantity depends on parameters $\g_i, i=1, 2, 3,$ and the relation between coordinates of the vector $\bn.$  We assume that $1/\al<\g_3<\g_1<\g_2<1$ and $n_3=\sv{n_1^\t n_2^\s}, $  where  $\t, \s>0.$ Thus we must investigate the quantity
\begin{equation}\label{Mn1n2}
M_{n_1, n_2}:=\max \left (n_1^{1-\g_1}n_2^{1-\g_2}, n_1^{\t(1-\g_3)}n_2^{\s(1-\g_3)}\right ).
\end{equation}
Let us define
\begin{equation}\label{tau0sigma0}
\t_0=\frac{1-\g_1}{1-\g_3}, \  \s_0=\frac{1-\g_2}{1-\g_3}, \ 0<\s_0<\t_0<1,
\end{equation}
and four sets ${\bar A}_j, 1\le j\le 4$ in the first quadrant of the $(\t, \s)$ plane:
$$
{\bar A}_1=\{\t\ge\t_0, \ \s\ge\s_0 \}, \quad {\bar A}_2=\{0<\t\le\t_0, 0<\s\le\s_0 \},
$$
$$
 {\bar A}_3=\{0<\t<\t_0, \ \s>\s_0 \}, \quad {\bar A}_4=\{\t>\t_0, \ 0<\s<\s_0 \}.
$$
It is easy to see that in the sets ${\bar A}_1$ and ${\bar A}_2$ one term from two terms in (\ref{Mn1n2}) is prevailing and we have
\begin{equation}\label{Mn1n21}
M_{n_1, n_2}=
\begin{cases}
n_1^{\t(1-\g_1)}n_2^{\s(1-\g_2)}, \text{if} \ (\t, \s)\in {\bar A}_1, \\
n_1^{1-\g_1}n_2^{1-\g_2}, \text{if} \ (\t, \s)\in {\bar A}_2.
\end{cases}
\end{equation}
In the sets ${\bar A}_3$ and ${\bar A}_4$ the situation is different and the maximum in (\ref{Mn1n2}) depends on the way how $(n_1, n_2)\to \infty$. We assume that $n_1=n, \ n_2=\sv{n^\r}, n_3=\sv{n^{\t+\r\s}}$ and let us consider $(\t, \s)\in {\bar A}_3$. Then we must find the quantity
\begin{equation}\label{MnA3}
M_{n}:=\max \left (n^{1-\g_1+\r(1-\g_2)}, n^{\t(1-\g_3)+\s(1-\g_3)}\right ).
\end{equation}
If we denote
\begin{equation}\label{rho0}
\r_0=\frac{\t_0-\t}{\s-\s_0},
\end{equation}
then it is easy to get that for $(\t, \s)\in {\bar A}_3$
\begin{equation}\label{MnA31}
M_{n}=
\begin{cases}
n^{1-\g_1+\r(1-\g_2)}, \text{if} \ 0<\r<\r_0, \\
n^{(\t+\r\s)(1-\g_3)}, \text{if} \ \r>\r_0.
\end{cases}
\end{equation}
Considering $(\t, \s)\in {\bar A}_4$,  assuming the same assumption $n_1=n, \ n_2=\sv{n^\r}, n_3=\sv{n^{\t+\r\s}}$, and using the same notation (\ref{rho0}) (note that in both sets ${\bar A}_3$ and ${\bar A}_4$ $0<\r_0<\infty$) we get
\begin{equation}\label{MnA4}
M_{n}=
\begin{cases}
n^{(\t+\r\s)(1-\g_3)}, \text{if} \ 0<\r<\r_0, \\
n^{1-\g_1+\r(1-\g_2)}, \text{if} \ \r>\r_0.
\end{cases}
\end{equation}
Thus, we got quite complicated behavior of the quantity (\ref{Mbn}). Assuming that $n_3=n_1^\t n_2^\s, $ we got four sets of parameters $\t, \s$, and in sets ${\bar A}_1$ and ${\bar A}_2$ the growth of coordinates $n_1$ and $n_2$ can be arbitrary and the quantity (\ref{Mbn}), which in this case becomes $M_{n_1, n_2}$, is given in (\ref{Mn1n21}). In sets ${\bar A}_3$ and ${\bar A}_4$ the growth of coordinates $n_1$ and $n_2$ cannot be arbitrary and we must assume $n_1=n, \ n_2=n^\r, n_3=n^{\t+\r\s}$, introducing new parameter $\r$. Then in each set ${\bar A}_3$ and ${\bar A}_4$ we got "boundary"value $\r_0$, the quantity (\ref{Mbn}) becomes $M_{n}$ and is given in (\ref{MnA31}) and (\ref{MnA4}). Having the expressions of the quantity (\ref{Mbn}) and remembering that $A_\bn=(n_1n_2n_3)^{1/\al}M_\bn$ (see (\ref{generalJn})) we can write down the normalization sequence $A_\bn$. Here it is appropriate to note that the value $\r_0=\r_0(\t, \s)$ is a function of $\t, \s$ and its behavior on the  sets ${\bar A}_3$ and ${\bar A}_4$ has the following properties. Let us take ${\bar A}_3$, then $\lim_{\t\to \t_0}\r_0(\t, \s)=0 $ for any fixed $\s$ and $\lim_{\s\to \s_0}\r_0(\t, \s)=\infty $ for any fixed $\t$, similar relations can be written for ${\bar A}_4.$

To complete the analysis of this example it would be necessary to find the limit distribution, which is quite complicated, but it is obtained in a standard way, therefore we mention the main step in finding the limit distribution. We can rewrite (\ref{fbnex8}) as follows
\begin{align}\label{fbnex8a}
M_\bn^{-1}f_{\bn}\left(\bu,\bt^{(j)}\right) &=\frac{z_{\g_1,n_1}z_{\g_2,n_2}}{M_\bn}\frac{U_{t_1^{(j)},\g_1}(\sv{n_1u_1},n_1)U_{t_2^{(j)},\g_2}(\sv{n_2u_2},n_2)}{z_{\g_1,n_1}z_{\g_2,n_2}}\ind{[0\leq\sv{n_3 u_3}  \leq n_3 t_3^{(j)}]}\\ \nonumber
                                  &+\frac{z_{\g_3,n_3}}{M_\bn}\frac{U_{t_3^{(j)},\g_3}(\sv{n_3u_3},n_3)}{z_{\g_3,n_3}}\ind{[0\leq\sv{n_i u_i}  \leq n_i t_i^{(j)}, \ i=1, 2]}
\end{align}
Having this expression we apply Proposition \ref{prop1}, find limits of ratios
$$
\frac{z_{\g_1,n_1}z_{\g_2,n_2}}{M_\bn},  \quad \frac{z_{\g_3,n_3}}{M_\bn}
$$
in sets ${\bar A}_i, 1\le i\le 4$, using expression of $M_\bn$ and our assumptions about relations between coordinates of $\bn$.

At the end of this example let us note that we had considered only one possible variant $1/\al<\g_3<\g_1<\g_2<1$. Clearly, it is possible to consider different location of parameters $\g_i$, and some variants will exhibit the scaling transition, while others will not. For example, changing only the location of $\g_3$ with respect to $\g_i, \ i=1, 2,$ we shall change only the location of the point $(\t_0, \s_0)$, namely, in the cases $1/\al<\g_1<\g_3<\g_2<1$ and $1/\al<\g_1<\g_2<\g_3<1$ we get $0<\s_0<1<\t_0$ and $1<\s_0<\t_0$, respectively. In the case $\g_i>1, i=1, 2, \ 1/\al<\g_3<1$ and $\sum_{i=0}^\infty a_j(i)\ne 0, \ j=1, 2,$ we have $M_\bn=z_{\g_3,n_3}, \ A_\bn=(n_1n_2)^{1/\al}n_3^{1/\al+1-\g_3}$ and there is no scaling transition. Using the terminology of \cite{Paul20}, one can say that in this case the r.f. under the consideration has zero memory in directions by the first two axis (horizontal plane) and positive memory in the vertical direction.

\section{Dependence structure of r.f. in the above examples}
The dependence structure of a linear r.f. is completely determined by the filter of the r.f. under consideration. Since the filters in the above presented examples are quite specific, the dependence structure in these examples is also  specific, and it can give some additional insight into phenomenon of the scaling transition.
 It is easy to  see that in Example 1  $X_{k,l}$ and $X_{m,n}$ are independent  if $k-m>0, l-n>0$ or $k-m<0, \ l-n<0$, in other cases these two values of the r.f. are dependent.
% ; in other words, $X_{k,l}$ depends only of the values of the r.f. on the vertical and horizontal lines going through the point $(k, l)$.
 Therefore, calculating covariances (in the case $\al=2$) or spectral covariances ($\al <2$; see \cite{PaulDam2}, where this quantity and other measures of dependence for r.f. are calculated) of this r.f. we shall get that big part of these quantities will be  zero. Let us take $\al=2$ in Proposition \ref{prop3} and let us consider the covariances $\rho (n, m):=\EE X_{0,0}X_{n, m}, (n, m)\in \ZZ^2$. It is easy to calculate that, as $|n|, |m| \to \infty$,
\begin{eqnarray}\label{cov1}
\rho(n, 0) & \sim & C|n|^{1-2\g_1}, \ \text {as}\  |n| \to \infty, \   \rho(0, m)\sim C|m|^{1-2\g_2},  \ \text {as}\ |m| \to \infty, \\ \nonumber
 \rho (n, m)& = & a_1(|n|)a_2(m) \ \ \text {if} \ n<0,  m>0, \\ \nonumber
  \rho (n, m)& =& a_1(n)a_2(|m|) \ \ \text {if} \ n>0,  m<0,
\end{eqnarray}
and $\rho (n, m)=\rho (-n, -m)=0$ if $n>0,  m> 0$. Here and in what follows $C$ stands for the constants, not the same at different places, which may be dependent on parameters $\g_i, \ i=1, 2$. Since in this example $1/2<\g_i<1$, we have the following relations
 \begin{equation}\label{cov2}
\sum_{(n, m) \in \ZZ^2}\rho (n, m)=\infty, \    \sum_{ m  \in \ZZ} \rho(0, m)=\infty, \  \sum_{n \in \ZZ} \rho(n, 0)=\infty.
\end{equation}
Usually for general stationary random fields with mean zero and finite variance, using the same notation $\rho (n, m)$ for covariances the relation
$\sum_{(n, m) \in \ZZ^2}|\rho (n, m)|=\infty$
is taken as definition of long-range dependence for the random field under consideration, while the relation
$\sum_{(n, m) \in \ZZ^2}|\rho (n, m)|<\infty$
serves as definition of short-range dependence. Remembering definition of directional memory for r.f. in \cite{Paul20}, it is possible to define long-range and short-range directional dependencies. We say that a stationary r.f. $\{X_{k, l}, \ (k, l)\in \ZZ^2 \}$ is long-range or short-range dependent in the horizontal direction, if for each fixed $m\in \ZZ$, the series
$\sum_{n \in \ZZ}|\rho (n, m)|$ is divergent or convergent, respectively.
Similarly we define directional dependence in the vertical direction. More generally, we can define both sorts of dependence for any direction, defined by means of rational numbers. Let $q=k/l$ be a fixed rational number (positive or negative) which defines a direction by means of the line $y=qx, x\in \RR.$ For  fixed $a, b\in \ZZ$, let us denote by $\cl(q,a, b)$ the set  $\{(lm+b, km+a): m\in \ZZ \}\subset \ZZ^2$ (integers $a, b$ are needed to include horizontal and vertical lines).
 \begin{definition} We say that a mean zero stationary r.f. $\{X_{k, l}, \ (k, l)\in \ZZ^2 \}$  with finite variance is long-range or short-range dependent in direction, defined by a rational number $q$, if for any fixed $a, b\in  \ZZ$ the series
$$
\sum_{(n, m)\in \cl(q, a, b)}|\rho (n, m)|
$$
is divergent or convergent, respectively. We say that this r.f. is long-range or short-range dependent if the series
$$
\sum_{(n, m)\in \ZZ^2}|\rho (n, m)|
$$
is divergent or convergent, respectively.
\end{definition}
Let us note that in \cite{Pilipauskaite} (see Remark 6.1 therein) vertical and horizontal long-range dependence was defined. Also it is necessary to note that sometimes in the definition of short-range dependence additionally it is required that the sum of covariances is not zero. If the sum of covariances is  zero, then we say that we have negative dependence. In our definition negative dependence is part of short-range dependence.

 Using these definitions we can say that (\ref{cov2}) means that the random field from Proposition \ref{prop3}  is long-range dependent and is  long-range dependent in both -vertical and horizontal - directions, but is short-range dependent in any other direction. Here it is necessary to note that it is easy to produce examples of a linear r.f. which has long-range dependence, but is short-range dependent in one  or even in both directions along axis, one such example will be r.f. from Example 3. Since in the case $\al=2$ the variance of a sum $S_{n, m}$ is expressed via covariances, analyzing the normalization constant $A_{n, m}$ it is easy to see that its growth depends on the way how $n, m$ tend to infinity.

Dependence structure of r.f. in Examples 2 and 3 is more complicated. Let as take a r.f. from Example 3, which has two points of scaling transition. Using the expression of the filter (\ref{cij2}) we have
$$
X_{k, l}=\sum_{i=1}^\infty a_i \e_{k-i, l} + \sum_{j=1}^\infty c_j \e_{k-j, l-j}.
$$
Therefore, it is easy to see that the random variable
$$
X_{0, 0}=\sum_{i=1}^\infty a_i \e_{-i, 0} + \sum_{j=1}^\infty c_j \e_{-j, -j}
$$
is independent only with variables $X_{k, l}$ if $k>l>0$ or $k<l<0$. For other combinations of indices $k, l$ it is not difficult to get the following relations for covariances:
%$$
%X_{0, l}=\sum_{i=1}^\infty a_i \e_{-i, l} + \sum_{j=1}^\infty c_j \e_{-j, l-j},\ X_{k, 0}=\sum_{i=1}^\infty a_i \e_{k-i, 0} + \sum_{j=1}^\infty c_j \e_{k-j, -j},\ X_{k, k}=\sum_{i=1}^\infty a_i \e_{k-i, k} + \sum_{j=1}^\infty c_j \e_{k-j, k-j},
%$$
%for all $l\ge 1, \ k\ge 1.$ Simple calculations give us
$$
\rho(0, l)=a_{|l|}c_{|l|}=(1+|l|)^{-\g_1-\g_2}, \    \rho(k, 0)=\sum_{i=1}^\infty a_i a_{|k|+i}\sim C|k|^{1-2\g_1}, \    \rho(k, k)=\sum_{i=1}^\infty c_i c_{|k|+i}\sim C|k|^{1-2\g_2}.
$$
Since $\g_1 +\g_2>1, \ -1,1-2\g_i <0, \ i=1, 2,$ we have
\begin{equation}\label{covex3}
\sum_{(n, m) \in \ZZ^2}\rho (n, m)=\infty, \    \sum_{ k  \in \ZZ} \rho(k, 0)=\infty, \ \sum_{ k  \in \ZZ} \rho(k, k)=\infty, \  \sum_{l \in \ZZ} \rho( 0, l)<\infty.
\end{equation}
These relations show that we have the example of a r.f. which we had mentioned above: it is long-range dependent but in vertical direction it is short-range dependent. In this example we have two directions - horizontal and diagonal -with long-range dependence. Analysis of normalizing constants $A_{n, m}$ shows more complicated behavior comparing with $A_{n, m}$ in Example 1.

Quite different situation is with Example 2. Since requiring the condition (\ref{sumZero}) we assume that $\sum_{i_1=0}^{\infty}\sum_{i_2=0}^{\infty}|c_{i_1,i_2}|<\infty$, it is not difficult to see that such linear random field is short-range dependent and, therefore, is short-range dependent in any direction. But it is easy to note that in this example we have the following relation:
\begin{equation}\label{covzero}
\sum_{(n, m) \in \ZZ^2}\rho (n, m)=0.
\end{equation}
Moreover, this relation is valid not only for our Example 2, but for general linear random field satisfying the condition (\ref{sumZero}). Namely, if a linear random field with absolutely summable filter $\{c_{i, j}, \ (i, j)\in \ZZ_+^2\}$ satisfies (\ref{sumZero}), the for such random field relation (\ref{covzero}) holds. The proof of this statement becomes very simple if we take filter defined on all $\ZZ^2$, then the proof follows from equalities
$$
\sum_{(n, m) \in \ZZ^2}\rho (n, m)=\sum_{(n, m) \in \ZZ^2}\sum_{(i, j) \in \ZZ^2}c_{i, j}c_{i+n, j+m}=
$$
$$
\sum_{(i, j) \in \ZZ^2}c_{i, j}^2  +\sum_{(i, j) \in \ZZ^2}c_{i, j} \sum_{(n, m) \in \ZZ^2, (n, m)\ne (0, 0)}c_{i, j}c_{i+n, j+m}=\left (\sum_{(i, j) \in \ZZ^2}c_{i, j} \right )^2.
$$
Returning  to the dependence structure in the Example 2 it is interesting to note, that despite of the relation (\ref{covzero}), sum of covariances over any line going trough the origin is positive. For example, denoting $A=\sum_{i=1}^\infty a_i, B=\sum_{i=1}^\infty b_i$  and recalling that $c_{0, 0}=a_0+b_0=-(A+B)$, we easily get
\begin{eqnarray*}
\sum_{n \in \ZZ}\rho (n, 0) & = & c_{0, 0}^2+2c_{0, 0}A+\sum_{i=1}^\infty a_i^2+2\sum_{i=1}^\infty a_i\sum_{n=1}^\infty a_{i+n}\\
                              & = & c_{0, 0}^2+2c_{0, 0}A+(\sum_{i=1}^\infty a_i)^2=(c_{0, 0}+A)^2>0
\end{eqnarray*}
In a similar way we can prove that $\sum_{(m) \in \ZZ}\rho (0, m)=(c_{0, 0}+B)^2>0$ and the sums $\sum_{(n) \in \ZZ}\rho (n, -n)$ and $\sum_{(n) \in \ZZ}\rho (-n, n)$ are also positive. That all these sums are positive can be explained by fact that in all these sums there is a big positive member $\rho_{0, 0}=\sum_{(i, j) \in \ZZ^2_+}c_{i, j}^2$. Most probably (but we did not verify) the same fact holds not only for our Example 2, but in the case of general filter if for all $(i, j)\ne (0, 0) \ c_{i, j}\ge 0$ and $c_{0, 0}=-\sum_{(i, j) \ne (0, 0)}c_{i, j}$.

Now let us look at the dependence structure in the same examples in the case $\al<2$. Note that the notions of long and short-range dependencies in the literature were used mainly for stationary r.f. with finite variance. This is due to the fact that our knowledge about dependence  for stable r.f. is quite limited. For long time there were only some results concerning dependence for stable r.f., see \cite{Samorod}, chapter 8.7, where codifference was calculated for some Takenaka r.f. defined on $\RR^d$. But the  dependence was measured in the following way: at first  r.f. was projected to a line, then for the obtained stable process on a line codifference, as a measure of dependence, was calculated.   Even for linear stable r.f. usual codifference (i.e., dependence between values of a r.f. at two points) was not investigated (reasons for that are explained in \cite{PaulDam2}). It turned out that the so-called $\al$-spectral covariance, introduced in  \cite{PaulDam2}, can serve as a measure of dependence and is quite successful substitute for the usual covariance in defining long-range and short-range dependencies. For definition of $\al$-spectral covariance we refer to \cite{PaulDam2}, here we recall only  that for  a linear r.f. (\ref{field}) $\al$-spectral covariance is given by formula
\begin{equation}\label{alphaspcor}
\rho_{\al}(n,m):=\rho_{\al}(X_{0,0}, X_{n,m})=\sum_{i=0}^\infty\sum_{j=0}^\infty c_{i, j}^{\langle \al/2\rangle}c_{i+n, j+m}^{\langle \al/2\rangle}, \ n>0, m>0,
\end{equation}
where $ x^{\langle a\rangle}=\abs{x}^a{\rm sign}(x)$.
Therefore, we suggest to classify  random fields (\ref{field}) in the same way as we classified  r.f. with finite variance. Namely, we say that a r.f. (\ref{field}) is long-range or short-range dependent (with respect to $\al$-spectral covariance) in direction, defined by a rational number $q$, if for any fixed $a, b\in  \ZZ$ the series
$$
\sum_{(n, m)\in \cl(q, a, b)}|\rho_{\al}(n, m)|
$$
is divergent or convergent, respectively. We say that this r.f. is long-range or short-range dependent (with respect to $\al$-spectral covariance), if  the series
$$
\sum_{(n, m)\in \ZZ^2}|\rho_{\al}(n, m)|
$$
is divergent or convergent, respectively. It is easy to see that this classification can be applied to r.f. indexed by $\ZZ^d, \ d\ge 2$ and  to any stationary r.f. for which we can define $\al$-spectral covariance.

Having expression (\ref{alphaspcor}), it is not difficult to verify, that the dependence structure in the examples of  Section 3 (except Example 2, since we do not know if some analogue of (\ref{covzero}) holds for $\al$-spectral covariance) in the case $\al<2$ remains the same as in the case $\al=2$. For example, taking the filter (\ref{cij}) from Example 1 and substituting this expression into (\ref{alphaspcor}), it is not difficult to get the following relations
$$
\rho_{\al}(0,m)=\sum_{i=1}^\infty a_2^{\al/2}(i) a_2^{\al/2}(m+i)\sim Cm^{1-2\b_2},   \  \rho_{\al}(n,0)=\sum_{i=1}^\infty a_1^{\al/2}(i) a_1^{\al/2}(n+i) \sim Cn^{1-2\b_1},
$$
where $\b_i:=\al \g_i/2$. Since $1/2<\b_i<\al/2<1$, we have
$$
\sum_{(n, m) \in \ZZ^2}\rho_{\al}(n, m)=\infty, \    \sum_{ m  \in \ZZ}\rho_{\al}(0, m)=\infty, \  \sum_{n \in \ZZ}\rho_{\al}(n, 0)=\infty,
$$
that is, exactly the same relations as in (\ref{cov2}), only with $\rho_{\al}(n, m)$ instead of $\rho (n, m)$.

Till now we had considered the dependence structure in examples of r.f. in the case $d=2$. Long-range and short-range dependence of stationary r.f. in the case $d=3$ can be defined in the same way as in the case $d=2$, but passing to the case $d=3$ we have more possibilities - we can sum covariances over points on lines or planes. For example, considering Example 4 and denoting $\rho (n_1, n_2, n_3):=\EE X_{0,0, 0}X_{n_1, n_2, n_3}, (n_1, n_2, n_3)\in \ZZ^3$ and taking $1/2<\g_1<\g_2<\g_3<1$, it is easy to see that non-zero covariances are on axes:
$$
\rho (n_1, 0, 0)\sim  C|n_1|^{1-2\g_1}, \quad \rho (0, n_2, 0)\sim  C|n_2|^{1-2\g_2}, \quad \rho (0, 0, n_3)\sim  C|n_3|^{1-2\g_3}.
$$
Thus, this r.f. is long-range dependent along each axis, but is short-range dependent along any line, going through origin and not coinciding with any of axes. Also it is easy to see that
$$
\sum_{(n_1, n_2) \in \ZZ^2}\rho (n_1, n_2, 0)=\infty, \quad \sum_{(n_1, n_3) \in \ZZ^2}\rho (n_1, 0, n_3)=\infty, \quad \sum_{(0, n_2, n_3) \in \ZZ^2}\rho (0, n_2, n_3)=\infty,
$$
but if we take the sum of covariances over plane, not containing any axis, we get the finite value.

In Example 5 (the case $\al=2$) we have r.f. with negative  and short-range dependence (condition (\ref{sum3dZero})), therefore short-range dependence will be over any  direction or plane. But there are specific relations for sums of covariances over the coordinate axes or coordinate planes. Let us denote
$$
V_1=\sum_{n_1\in \ZZ}\rho (n_1, 0, 0), \quad V_2=\sum_{n_2\in \ZZ}\rho (0, n_2, 0), \quad V_3=\sum_{n_3\in \ZZ}\rho (0, 0, n_3),
$$
$$
U_{3}=\sum_{(n_1, n_2)\in \ZZ^2}\rho (n_1, n_2, 0), \quad U_{2}=\sum_{(n_1, n_3)\in \ZZ^2}\rho (n_1, 0, n_3), \quad U_{1}=\sum_{(0, n_2, n_3)\in \ZZ^2}\rho (0, n_2, n_3),
$$
$$
A_j=\sum_{i=1}^\infty a_j (i), \quad B_j=\sum_{i=1}^\infty a_j^2 (i), \ j=1, 2, 3, \quad \rho_0=\rho (0, 0, 0), \ {\bar c}_0=c_{(0, 0, 0)}.
$$
We have the following result.
\begin{prop}\label{prop7} In Example 5 in the case of finite variance we have the following relations:
\begin{equation}\label{covzero0}
\sum_{(n_1, n_2, n_3)\in \ZZ^3}\rho(n_1, n_2, n_3)=0
\end{equation}
\begin{equation}\label{covzero1}
V_k=\left (\sum_{j=1, j\ne k}^3 A_j \right )^2 + \left (\sum_{j=1, j\ne k}^3 B_j \right )>0, \ k=1, 2, 3,
\end{equation}
\begin{equation}\label{covzero2}
U_k= A_k^2 +B_k>0,
\end{equation}
i.e., all sums of covariances over the coordinate axes or coordinate planes are positive and only sum of covariances over $\ZZ^3$ is equal to zero.
\end{prop}

{\it Proof of Proposition \ref{prop7}}. We prove (\ref{covzero1}) for $k=1$, since other two relations can be proved in the same way. We have
$$
\rho (n_1, 0, 0)={\bar c}_0c_{(n_1, 0, 0)}+\sum_{i=1}^\infty a_1(i)c_{(i+n_1, 0, 0)}
$$
and
$$
V_1=\sum_{n_1\in \ZZ}\rho (n_1, 0, 0)=\rho_0+ \left (\sum_{n_i=-\infty}^{-1} +\sum_{n_1=1}^\infty \right )\rho (n_1, 0, 0).
$$
From definition of covariance we  have $\rho_0={\bar c}_0^2+B_1+B_2+B_3$ and easy calculations give us
$$
\left (\sum_{n_i=-\infty}^{-1} +\sum_{n_1=1}^\infty \right )\rho (n_1, 0, 0)=2\left ({\bar c}_0 A_1+\sum_{i=1}^\infty\sum_{m=1}^\infty a_1(i) a_1(i+m)\right ).
$$
Therefore, we have
\begin{eqnarray*}
V_1 &=& {\bar c}_0^2+2{\bar c}_0 A_1+\sum_{i=1}^\infty a_1^2(i)+2\sum_{i=1}^\infty\sum_{m=1}^\infty a_1(i) a_1(i+m)+B_2+B_3\\
    &=& ({\bar c}_0+A_1)^2+B_2+B_3.
\end{eqnarray*}
Taking into account the relation ${\bar c}_0=-(A_1+A_2+A_3)$ we get  (\ref{covzero1}) with $k=1.$

Now we prove (\ref{covzero2}) with $k=3$ (the proof for $k=1, 2$ is similar). Let us note that $\rho (n_1, n_2, 0)=0$, if $n_1n_2>0$ and $\rho (n_1, n_2, 0)=a_1(|n_1|)a_2(|n_2|)$, if $n_1n_2<0$.
Therefore, we can write
$$
U_3=V_1+V_2-\rho_0 +\left (\sum_{n_1=-\infty}^{-1}\sum_{n_2=1}^\infty +\sum_{n_1=1}^\infty\sum_{n_2=-\infty}^{-1} \right )\rho (n_1, n_2, 0).
$$
It is easy to see that
$$
\left (\sum_{n_1=-\infty}^{-1}\sum_{n_2=1}^\infty +\sum_{n_1=1}^\infty\sum_{n_2=-\infty}^{-1} \right )\rho (n_1, n_2, 0)=A_1A_2,
$$
therefore, using expressions (\ref{covzero1}) for $V_1, V_2$ and the relation $\rho_0={\bar c}_0^2+B_1+B_2+B_3$ we get
$$
U_3=(A_2+A_3)^2 + (A_1+A_3)^2+B_1+B_2+2B_3-{\bar c}_0^2-(B_1+B_2+B_3)+2A_1A_2.
$$
From this relation we easily get (\ref{covzero2}) with $k=3$.

Although we had proved (\ref{covzero0}) earlier, writing the identity
$$
\sum_{(n_1, n_2, n_3)\in \ZZ^3}\rho(n_1, n_2, n_3)=U_1+U_2+U_3-V_1-V_2-V_3+\rho_0
$$
and using (\ref{covzero1}) and (\ref{covzero2}), we can verify (\ref{covzero0}).
\halmos

%\vfill
%\eject
%%
%
%
%

%\halmos
%\bibliographystyle{plain}
%\bibliography{probability13}

\begin{thebibliography}{10}

\bibitem{PaulDam2}
J. Damarackas and V. Paulauskas.
Spectral covariance and limit theorems  for random fields  with infinite variance.
{\em J. Multivar. Analysis}, 153: 156--175, 2017.



\bibitem{DavPau1}
Yu. Davydov and V. Paulauskas.
Lamperti type theorems for random fields.
{\em Teor. verojatn. i primenen.},  63: 520--544, 2018.

\bibitem{Gnedenko}
B.V. Gnedenko and A.N. Kolmogorov.
{\em Limit distributions for sums of independent random variables}.
Addison-Wesley, Reading, MA, 1968


\bibitem{Lahiri}
S.N. Lahiri and P.M. Robinson.
 Central limit theorems for long range dependent spatial linear processes.
{\em Bernoulli}, 22: 345--375, 2016.


\bibitem{Lamperti}
 J. Lamperti.
Semi-stable stochastic processes.
{\em Trans American Math. Soc.}, 104: 62--78, 1962.



\bibitem{Paul20}
V. Paulauskas.
Some remarks on definition of memory for stationary random processes and fields.
{\em Lith. Math. J.}, 56: 229--250, 2016.



\bibitem{Pilipauskaite}
V. Pilipauskait{\. e} and D. Surgailis.
Scaling transition for non-linear  random fields with long range dependence.
 {\em Stochastic Process. Appl.}, 127: 2751--2779, http://dx.doi.org/10.1016/j.spa.2016.12.011, 2017



\bibitem{Puplinskaite2}
D.~Puplinskait{\. e} and D.~Surgailis.
 Scaling transition for long-range dependent {Gaussian} random fields.
 {\em Stochastic Process. Appl.}, 125:2256--2271, 2015.

\bibitem{Puplinskaite1}
D.~Puplinskait{\. e} and D.~Surgailis.
 Aggregation of autoregressive random fields and anisotropic
  long-range dependence.
 {\em Bernoulli}, 22:2401--2441, 2016.


\bibitem{Samorod}
G.~Samorodnitsky and M.~Taqqu.
 {\em Stable non-Gaussian Random Processes. Models with Infinite
  Variance}.
 Chapman{\&} Hall, New York, 1994.

\bibitem{Surg}
D. Surgailis.
Anisotropic scaling limits of long-range dependent random fields on ${Z}^3$.
{\em J. Math. Anal. Appl.},  472: 328--351, 2019.


\end{thebibliography}

\end{document}